\documentclass[10pt,draftcls,onecolumn]{IEEEtran}
\usepackage{graphicx}
\usepackage[ruled,vlined]{algorithm2e}
\usepackage{url}
 \usepackage[pdftex]{color}
\usepackage{amsmath}
\usepackage{amsfonts} 
\usepackage{array}
\usepackage{cancel}

\newtheorem{propo}{Proposition}[section]

\newtheorem{rem}{Remark}
\usepackage{soul}
\usepackage{tabularx}
\usepackage{multirow}

\newcommand{\beq}{\begin{equation}}
\newcommand{\eeq}{\end{equation}}
\newcommand{\bea}{\begin{eqnarray}}
\newcommand{\eea}{\end{eqnarray}}

\newcommand{\bsea}{\begin{subeqnarray}}
\newcommand{\esea}{\end{subeqnarray}}
\newcommand{\nn}{\nonumber}

\newcommand{\qed}{\hfill $\Box$ \vskip 2ex}

\def\bmat{\left[ \begin{array}}
\def\emat{\end{array} \right]}
\DeclareMathOperator{\tr}{tr} 

\definecolor{Royalblue}{cmyk}{1,0.30,0.2,0.2}

\newcommand{\alg}[1]{\begin{align}#1\end{align}}

\newcommand{\Rs}{\mathbb{R}} 
\newcommand{\Es}{\mathbb{E}} 
\newcommand{\Nc}{\mathcal{N}} 
\newcommand{\Bc}{\mathcal{B}}

\newcommand{\Vc}{\mathcal{V}}
\newcommand{\Sc}{\mathcal{S}}
\newcommand{\Cc}{\mathcal{C}}
\newcommand{\Ac}{\mathcal{A}}
\newcommand{\Fc}{\mathcal{F}}
\newcommand{\Gc}{\mathcal{G}}
\newcommand{\Qc}{\mathcal{Q}}
\newcommand{\Hc}{\mathcal{H}}
\newcommand{\Rc}{\mathcal{R}}
\newcommand{\Pc}{\mathcal{P}}

\begin{document}

\title{Distributed Kalman Filtering\\ under Model Uncertainty}

\author{Mattia~Zorzi \thanks{M. Zorzi is with the Department of Information
Engineering, University of Padova, Padova, Italy; e-mail: {\tt\small zorzimat@dei.unipd.it} .}}

\markboth{DRAFT}{}

\maketitle

\begin{abstract}
We study the problem of distributed Kalman filtering for sensor networks in the presence of model uncertainty. More precisely, we assume that the actual state-space model belongs to a ball, in the Kullback-Leibler topology, about the nominal state-space model and whose radius reflects the mismatch modeling budget allowed for each time step. We propose a distributed Kalman filter with diffusion step  which is robust with respect to the aforementioned model uncertainty. Moreover, we derive the corresponding least favorable performance. Finally, we check  the effectiveness of the proposed algorithm in the presence of uncertainty through a numerical example.    
 \end{abstract}


\section{Introduction}

Significant advances in science and technology have led to a large number of problems that involve numerous sensors, i.e.  a sensor network (SN), taking measurements and a state process which needs to be estimated from such measurements. Just to mention a few of these problems: area  surveillance, region monitoring, target tracking and electrical power grid analysis. These problems can be (in principle) solved by using Kalman filtering 
equipped with all the observations coming from the SN. On the other hand, such centralized strategy is impractical or impossible to implement. 
Indeed, it requires a large amount of energy for communications
among the central node, i.e. the one which computes an estimate of the state process, and the sensors. In order to overcome this difficulty, distributed strategies have gained rapidly increasing interest in the last few years, see for instance \cite{distr_cons_TCN18,distr_joint_TCN15,av_stateobs_TCN17,fully_distr_TCN18}. The latter represent an attractive alternative because they require fewer communications and allow parallel processing. The simpleminded distributed version of the Kalman filter assumes that each node can compute an estimate of the state by using only the observations coming from its neighbors. On the other hand, such an approach provides poor performances compared to the ones of the centralized approach. A remarkable improvement has been gained by the so called consensus-based distributed Kalman filters \cite{spanos2005,Saber_CDC2005,Saber_CDC2007,Carli2008}.
Such approaches require multiple communication iterations during each sampling time interval: for instance, in the first iteration the nodes exchange their observation and compute their local estimate; in the second iteration the nodes exchange their local estimates and construct the final estimate based on consensus schemes. A further improvement has been given by diffusion-based strategies, \cite{cattivelli2010,yang2016diffusion}, wherein the consensus law is replaced by a convex combination of the local estimates. It is worth noting that there also exist high performing distributed Kalman filters based on different principles. For instance it is possible to compute the estimate through a fusion center which merges the local estimates of the nodes, \cite{song2007optimal,xu2012optimal,SUN20041017,doi:10.1080/00207721.2010.502601}.

Kalman filtering is based on nominal state-space models. On the other hand, the latter are just an approximation of the underlying system, thus the resulting estimate could lead to poor performances in practice. To address this model uncertainty  issue many robust Kalman filtering strategies have been proposed. The most popular one is risk-sensitive filtering, \cite{boel2002robustness,ROBUSTNESS_HANSENSARGENT_2008}, wherein large errors are severely penalized according to the so called risk sensitivity parameter: the larger the latter is the more large errors are penalized. Here, we consider the robust Kalman filter introduced by Levy and Nikoukhah \cite{LEVY_NIKOUKHAH_2004,ROBUST_STATE_SPACE_LEVY_NIKOUKHAH_2013}. In such approach, the actual state space model belongs to a ball, according to the Kullback-Leibler divergence, about the nominal model and with radius, say tolerance, which represents the mismatch modeling budget allowed for each time step. 
Then, the optimal estimator is designed according to the least favorable model in this ball. It is turns out that the latter obey to a Kalman-like recursion. Furthermore, it can be understood as a generalization of risk-sensitive filtering: such a filter can be rewritten as a risk-sensitive filter with a time varying risk-sensitivity parameter. Finally, it worth noting that the aforementioned approach can be extended to a family of balls formed by using the $\tau$-divergence, \cite{OPTIMALITY_ZORZI,STATETAU_2017}.

It is then natural to wonder how to perform distributed Kalman filtering for SNs in the presence of model uncertainty. Many papers consider model uncertainty in terms of missing observations. Such a situation well describes communication problems among sensors: for instance, an $H_\infty$-consensus problem for SNs with multiple missing measurements has been considered in  \cite{SHEN20101682}. On the other hand, only few papers addressed the problem of model uncertainty in a broader sense to the best of the author's knowledge. For instance, a
 distributed Kalman filtering fusion strategy with random state transition and measurements matrices has been considered in \cite{Luo_2008}.

The contribution of the present is about distributed Kalman filtering for SNs wherein the model uncertainty is characterized by a ball, in the Kullback-Leibler divergence topology, about the nominal model. More precisely, we propose a distributed version of the robust Kalman filter introduced by Levy and Nikoukhah. The proposed algorithm is a robust version of the distributed Kalman filter with diffusion step in \cite{cattivelli2010}. Then, we derive the least favorable performance of these filters. Similarly to the centralized case, the latter can be characterized over a finite simulation horizon as follows: first, a forward recursion is required to compute the optimal gains of centralized Kalman filter; then, a backward recursion is required to compute the least favorable model; finally, a forward recursion is required to compute the performance of the distributed robust Kalman algorithm. We show that the average least favorable mean square deviation across the network converges to a finite constant value provided that: the tolerance is sufficiently small; reachability and local observability hold. 
 Moreover, we show that it is very likely that the proposed robust distributed filters perform better than the standard ones provided that the tolerance is sufficiently large.
Finally, we compare the proposed algorithms with the standard ones through a numerical example.

The outline of the paper is as follows. In Section \ref{sec_back} we introduce the background about robust Kalman filtering. In Section
\ref{sec_distr_K} we present the distributed robust Kalman filtering algorithm. In Section \ref{sec_lf} we analyze the least favorable performance of the proposed algorithms. In Section \ref{sec_num} we present a numerical example which compares the proposed algorithm with the standard ones. Finally, in Section \ref{sec_concl} we draw the conclusions.

\section{Background}\label{sec_back}

Consider the nominal state-space model
\alg{ \label{nomi_mod}x_{t+1}&=Ax_t+\Gamma_B u_t\nn\\
y_t&=C x_t+\Gamma_Du_t}
where $x_t\in\Rs^n$ is the state process, $y_t\in \Rs^{pN}$ is the observation process, $u_t$ is normalized white Gaussian noise (WGN), i.e. $\Es[u_tu_s^T]=I \delta_{t-s}$ where $\delta_t$ denotes the Kronecker delta function.
We assume that $u_t$ is independent of  the initial state $x_0$. The latter is Gaussian distributed with mean $\hat x_0$ and covariance matrix $V_0$. Model (\ref{nomi_mod})
is characterized by the nominal transition probability density of $z_t:=[\, x_{t+1}^T\; y_t^T\,]^T$ given $x_t$ which is denoted by $\phi_t(z_t|x_t)$. We assume that $\Gamma_B\Gamma_D^T=0$, i.e. the noise entering in the state process is independent of the noise entering in the observation process.   We assume that $u_t$ affects all the components of the dynamics and observations in (\ref{nomi_mod}). Such assumption is necessary whenever entropy-like indexes are used to measure the proximity of statistical models, as in our case, otherwise these indexes take infinite value. Accordingly, the matrix $[\,\Gamma_B^T \;\Gamma_D^T\,]^T$ is full row rank, and without loss of generality we can assume that $[\,\Gamma_B^T \;\Gamma_D^T\,]^T$ is a square (and thus invertible) matrix of dimension $pN+n$. Indeed, we can always compress the column space of such a matrix and remove the noise components which do not affect model in (\ref{nomi_mod}).
Accordingly, the state-space model in (\ref{nomi_mod}) is reachable; moreover we also assume it is observable.

Let $\tilde \phi_t(z_t|x_t)$ be the (unknown) actual transition probability density of $z_t$ given $x_t$. In order to account the fact that the nominal model does not coincide with the actual model, we assume that $\tilde \phi_t$ belongs to the closed ball about $\phi_t$:
\alg{\label{ball_B}\mathbb{B}_t:=\{\, \tilde \phi_t \hbox{ s.t. } \tilde\Es[\log (\tilde \phi_t / \phi_t)| Y_{t-1}]\leq c\,\}}  
where 
\alg{\label{def_ball_E}\tilde\Es&[\log (\tilde \phi_t / \phi_t)| Y_{t-1}]\nn\\ &:=\int \int \tilde \phi_t(z_t|x_t) \tilde f_t(x_t|Y_{t-1})\log\left(\frac{\tilde \phi_t(z_t|x_t)}{\phi_t(z_t|x_t)}\right) \mathrm{d}z_t \mathrm{d}x_t}
and $Y_{t-1}:=\{\,y_s,\;\; s=1\ldots t-1 \,\}$. The latter represents the relative entropy between the actual and the nominal transition densities $\tilde \phi_t(z_t|x_t)$ and $ \phi_t(z_t|x_t)$ at time $t$, respectively, and $\tilde f_t(x_t|Y_{t-1})\sim\Nc(\hat x_t, V_t)$ is the actual conditional probability density of $x_t$ given the past observations $Y_{t-1}$. Finally, parameter $c>0$ is called tolerance and represents the mismatch modeling budget allowed for each time step.

Given the nominal model in (\ref{nomi_mod}), a robust estimator of $x_{t+1}$ given $Y_{t}$ is obtained by solving the following mini-max problem
\alg{\label{minimax_pb} \hat x_{t+1} &= \underset{g_t\in\mathbb G_t}{\mathrm{argmin} \,}\underset{\tilde \phi_t\in\mathbb B_t}{\mathrm{max}} \tilde \Es[\| x_{t+1}-g_t\|^2|Y_{t-1}]} where \alg{ \label{minimax_pbE}\tilde \Es[\| & x_{t+1}-g_t\|^2|Y_{t-1}]\nn\\ &:=\int \int \|x_{t+1}-g_t\|^2\tilde \phi_t(z_t|x_t) \tilde f_t(x_t|Y_{t-1})\mathrm{d}x_t \mathrm{d}z_{t}}
represents the mean square error of the estimator $g_t$ which is a function of $y_t$ and $Y_{t-1}$. $\mathbb G_t$ denotes the set of all estimators $g_t$ such that $\tilde \Es[\|g_t\|^2]$ is finite for any $\tilde \phi_t \in\mathbb B_t$. Roughly speaking, such estimator is designed according to the least
favorable model whose mismodeling budget allowed is expressed at each time step. This way to characterize model uncertainty is better than expressing the uncertainty over the entire simulation interval. Indeed, in the latter case the maximizer has the possibility to identify the moment where the dynamic of model (\ref{nomi_mod}) is most susceptible  to distortion and to allocate most of the distortion budget specified by the tolerance to this single element of the model, that is a situation which is pretty unrealistic.

In \cite{ROBUST_STATE_SPACE_LEVY_NIKOUKHAH_2013} it has been proved that the estimator solving the mini-max problem (\ref{minimax_pb}) obeys the Kalman-like recursion:
\alg{ \label{rob_est} G_t&=A V_t C^T(C V_t C^T+\Gamma_D \Gamma_D^T)^{-1} \nn\\
\hat x_{t+1}&=A\hat x_t+G_t(y_t-C \hat x_t) \nn\\
 P_{t+1}&= A(V_t^{-1} +C^T (\Gamma_D \Gamma_D^T)^{-1}C)^{-1}A^T+\Gamma _B\Gamma_B^T\nn\\
 \hbox{Find } & \hbox{$\theta_t$ s.t. $\gamma(P_{t+1},\theta_t)=c$}\nn\\
V_{t+1}&=(P_{t+1}^{-1}-\theta_t I)^{-1}
     }
where
\alg{\label{def_gamma}\gamma(P,\theta):=\log\det(I-\theta P)+\tr((I-\theta P)^{-1}-I).}
The so called risk sensitivity parameter, \cite{LEVY_ZORZI_RISK_CONTRACTION}, $\theta_t>0$ does always exist and it is unique given $P_{t+1}$ and $c$, moreover it can be computed efficiently by using a bisection algorithm. In the limit case $c=0$, i.e. there is no uncertainty, then $\theta_t=0$ and (\ref{rob_est}) becomes the usual Kalman filter.

\begin{rem} The robust filtering paradigm in (\ref{minimax_pb}) can be extended to nominal state-space models with time-varying parameters and tolerance. On the other hand, to ease the introduction of the corresponding distributed algorithms we stick to the constant parameters and tolerance case. 
\end{rem}

\section{Distributed Robust Kalman filter} \label{sec_distr_K}

Consider a network of $N$ nodes and in each node there is one sensor. We say that two nodes are connected if the corresponding sensors can communicate directly with each other. A node is always connected with itself. The neighborhood of node $k$, i.e. the set of nodes connected with $k$, is denoted by $\Nc_k$, in particular $k\in\Nc_k$. The corresponding $N\times N$ adjacency matrix $J=[j_{lk}]_{lk}$ is defined as
\alg{j_{lk}:=\left\{\begin{array}{ll}  1,& \hbox{if } l\in\Nc_k \\ 0,& \hbox{otherwise.}\end{array}\right.} The number of neighbors of node $k$ is denoted by $n_k$. Every node at time $t$ collects a measurement $y_{k,t}\in \Rs^{p}$ whose underlying model is unknown. The nominal model takes the form:
\alg{ \label{distr_mod}x_{t+1}&=Ax_t+B w_t\nn\\
y_{k,t}&=C_k x_t+D_kv_{k,t} \; \; k=1\ldots N}
where $w_t$ and $v_{k,t}$ $k=1\ldots N$ are independent WGNs such that $\Es[w_tw_t^T]=I$, $\Es[w_tv_{k,t}^T]=0$ and
$\Es[v_{l,t}v_{k,s}^T]=I\delta_{k-l}\delta_{t-s}$. It is worth noting that (\ref{distr_mod}) can be rewritten as (\ref{nomi_mod}) with $y_t=[\,y_{1,t}^T \ldots \, y_{N,t}^T\,]^T$, $u_t=[\,w_t^T \,v_{t}^T \,]^T$, $v_t=[\,v_{1,t}^T \ldots \, v_{N,t}^T\,]^T$, $\Gamma_B=[\,B \; 0\,]$, $\Gamma_D=[\,0 \; D\,]$,
\alg{C=\left[\begin{array}{c}C_1  \\ \vdots \\ C_N\end{array}\right], \; \;
D=\mathrm{diag}(\, D_1,\, \ldots \, ,  D_N)} 
where $\mathrm{diag}$ is the linear operator which constructs a block diagonal matrix whose blocks are the ones specified in the argument. We also define $R:=DD^T$, $R_l:=D_lD_l^T$ with $l=1\ldots N$, and 
\alg{S_{tot}:=C^TR^{-1}C=\sum_{l=1}^N C_l^TR_l^{-1}C_l.}
Accordingly, the filtering gain for model (\ref{distr_mod}) can be written as
\alg{G_t&=AV_tC^T (C V_t C^T+R)^{-1}= AV_tC^T\times \nn\\ & \hspace{0.4cm}\times (R^{-1}- R^{-1}C (V_t^{-1}+C^TR^{-1}C )^{-1}  C^T R^{-1})\nn\\
& = A (V_t^{-1}+C^TR^{-1}C )^{-1}  C^T R^{-1}\nn\\
&= A \left(V_t^{-1}+S_{tot} \right)^{-1}  C^T R^{-1}\nn}
and thus
\alg{G_ty_t&=A \left(V_t^{-1}+S_{tot} \right)^{-1}  \sum_{l=1}^N C_l^TR_l^{-1}y_{l,t} \nn\\
G_tC&= \left(V_t^{-1}+S_{tot}\right)^{-1}  S_{tot}.\nn}

In distributed Kalman filtering under model uncertainty, the aim is to compute for every node $k$ a prediction of the state $x_t$ while sharing the data only with its neighbors $l\in \Nc_k$ and taking into account that (\ref{distr_mod}) does not coincide with the actual model. In what follows, the one step-ahead prediction of $x_t$ at node $k$ is denoted by $\hat x_{k,t}$. It is not difficult to see that the robust Kalman filter for model (\ref{distr_mod}), i.e. the node $k$ has access to all measurements across all the nodes in the network, can be written as
\alg{ \label{rob:est2} 
\hat x_{k,t+1}&=A\hat x_{k,t}+A(V_{k,t}^{-1}+S_{tot} )^{-1}  \times \nn\\ & \hspace{2.2cm}\times\sum_{l=1}^N C_l^TR_l^{-1}(y_{l,t}-C_l\hat x_{k,t})\nn\\
 P_{k,t+1}&= A\left(V_{k,t}^{-1}+S_{tot}\right)^{-1} A^T+BB^T\nn\\
 \hbox{Find } & \hbox{$\theta_{k,t}$ s.t. $\gamma(P_{k,t+1},\theta_{k,t})=c$}\nn\\
V_{k,t+1}&=(P_{k,t+1}^{-1}-\theta_{k,t} I)^{-1} } where $\hat x_{k,t}=\hat x_t$ and $V_{k,t}=V_t$ with $k=1\ldots N$.

Therefore in the case that  the node $k$ has not access to all measurements across all nodes in the network, we would obtain a state prediction $\hat x_{k,t}$ of $x_t$ which is as close as to the global state prediction. 

 \subsection{Robust Kalman filter with diffusion step}
We assume that a node $k$ has access to the measurements  of its neighbors $\Nc_k$. The corresponding nominal state-space model is 
\alg{x_{t+1}&= Ax_t+Bw_t\nn\\
y_{l,t} &= C_lx_t+D_l v_{l,t},\; \; l\in\Nc_k.} 
The latter can be rewritten as a state-space model $(A,\Gamma_B,C_k^{loc},\Gamma_{D_k^{loc}})$ with 
input noise $u_{k,t}^{loc}=[\,w_t^T\; \,(v_{k,t}^{loc})^T\,]^T$ and output $y_{k,t}^{loc}$ where $v_{k,t}^{loc}$ and $y_{k,t}^{loc}$
are obtained by stacking $v_{l,t}$ and $y_{l,t}$, respectively, with $l\in\Nc_k$. Moreover, $\Gamma_B=[\,B \;0\,]$, $C_k^{loc}$ is obtained by stacking $C_l$ with $l\in\Nc_k$, $\Gamma_{D_k^{loc}}=[\,0 \;D_k^{loc}\,]$
 and $D_k^{loc}$ is a block diagonal matrix whose main blocks are $D_l$ with $l\in\Nc_k$. We also define $R_k^{loc}:=D_k^{loc}(D_k^{loc})^T$, $S_k:=(C_k^{loc})^T(R_k^{loc})^{-1}C_k^{loc}$ and thus
\alg{S_k=\sum_{l\in\Nc_k}C_l^TR_l^{-1} C_l.}
Accordingly, the one-step ahead predictor $\hat x_{k,t}$ of $x_t$ at node $k$ is given by (\ref{rob:est2}) where the terms for which $l\notin\mathcal N_k$ are discarded. Then,  the local prediction
$\hat x_{k,t+1}$  can be understood as an intermediate
local prediction of $x_t$ at node $k$. In what follows we denote such intermediate prediction as $\psi_{k,t+1}$. Then, the idea is to update the prediction at node $k$ not only in terms of $\psi_{k,t+1}$, but also in terms of $\psi_{l,t+1}$ with $l\in\Nc_k$. More precisely, we consider a matrix $W=[w_{lk}]_{lk}\in\Rs^{N\times N}$ such that 
\alg{\label{prop_W} &w_{lk}\geq 0 \hbox{ and } w_{lk}=0 \hbox{ if } l\notin \Nc_k\nn\\
& \sum_{l\in\Nc_k} w_{lk}=1. }
Then, the final prediction at node $k$ is obtained by the so called diffusion step, \cite{cattivelli2010}:
\alg{\label{diff_step}\hat x_{k,t+1}=\sum_{l\in\Nc_k} w_{lk}\psi_{l,t+1}.}
Therefore, in the diffusion algorithm a node $k$ exploits the information of the neighbors in terms of $y_{l,t}$ and $\psi_{l,t+1}$. The resulting procedure is outlined in Algorithm \ref{rob_est_diff}.
\begin{algorithm}\label{rob_est_diff}
 \small
 \SetKwInOut{Input}{Input}\SetKwInOut{Output}{Output}
 \Input{$\hat x_{k,t}$, $V_{k,t}$, $y_{k,t}$ with $k=1\ldots N$}\smallskip
 \Output{$\hat x_{k,t+1}$, $V_{k,t+1}$ with $k=1\ldots N$}\smallskip
 \textbf{Incremental step.} Compute at every node $k$:\\ \smallskip
  \hspace*{0.3cm} $ \psi_{k,t+1}=$\\
 \smallskip \hspace*{0.8cm} $A\hat x_{k,t}+ A(V_{k,t}^{-1}+S_k )^{-1}  \sum_{l\in\Nc_k} C_l^TR_l^{-1}(y_{l,t}-C_l\hat x_{k,t})$\\ \smallskip
  \hspace*{0.3cm} $P_{k,t+1}= A(V_{k,t}^{-1}+S_k)^{-1} A^T+BB^T$\\\smallskip
  \hspace*{0.3cm} Find $\theta_{k,t}$ s.t. $\gamma(P_{k,t+1},\theta_{k,t})=c$\\\smallskip
 \hspace*{0.3cm} $V_{k,t+1}=(P_{k,t+1}^{-1}-\theta_{k,t} I)^{-1} $\\ \smallskip
\textbf{Diffusion step.} Compute at every node $k$:\\ \smallskip
 \hspace*{0.3cm} $\hat x_{k,t+1}=\sum_{l\in\Nc_k} w_{lk} \psi_{l,t+1}$
   \caption{Distributed Robust Kalman filter with diffusion step at time $t$}
\end{algorithm} As explained in \cite{cattivelli2010}, the diffusion step (\ref{diff_step}) is motivated by the fact that the centralized prediction $\hat x_{t+1}$ can be approximated by a local convex combination of $\psi_{l,t+1}$. It is worth noting that in the case that $c=0$, i.e. there is no mismatch between the actual and the nominal model,  then $\theta_{k,t}=0$ for any $t$ and $k$ so that, we obtain the diffusion algorithm proposed in \cite{cattivelli2010}. In the case that 
\alg{\label{W_cons}w_{lk}:=\left\{\begin{array}{ll}\varepsilon,  & \hbox{if } l\neq k,\;l\in\Nc_k \\ 1-\varepsilon(n_k-1),  &  \hbox{if } l=k,\;l\in\Nc_k \\
0,  &  \hbox{otherwise} 
  \end{array}\right.}
we obtain a consensus-based update where $\varepsilon>0$ is the consensus parameter. Indeed, in the case that $c=0$ and $W$ is designed as in (\ref{W_cons}) we obtain the distributed consensus-based algorithm proposed in \cite{Saber_CDC2007}. Finally, in the case that $W=I$ we obtain the robust version of the local Kalman filter  \cite[p. 329]{KAILATH_OPTIMAL_FILTERING}.

It is worth noting that the mismatch modeling budget $c$ in Algorithm \ref{rob_est_diff} coincides with the one of the centralized filter. Such a choice does not guarantee that the least favorable model computed at node $k$ coincides with the one of the centralized filter. On the other hand, we will see that, under large deviations of the least favorable model of the centralized problem, it is very likely that the predictor at node $k$ using Algorithm \ref{rob_est_diff} performs better than the one in \cite{cattivelli2010},  see Section \ref{sec_opt} for more details.

\begin{rem} In some cases we may have a state-space model of the form
\alg{x_{t+1}&=Ax_t+\Gamma_B u_t+r_t\nn\\
y_t&=C x_t+\Gamma_Du_t\nn}
where $r_t$ is a deterministic process. In 
\cite{RS_MPC_IET} it was shown 
that the corresponding centralized robust Kalman filter still obeys the Kalman-like recursion (\ref{rob_est}) where the prediction update is replaced by 
\alg{\hat x_{t+1}&=A\hat x_t+G_t(y_t-C \hat x_t)+r_t.} 
Then, it is not difficult to see that the distributed algorithms presented in this section still hold in this case. The unique difference is that we need to add in the prediction update of each node the term $r_t$.
\end{rem}

\section{Least favorable performance} \label{sec_lf}
In this section we analyze the performance of the distributed algorithm with diffusion step introduced in Section \ref{sec_distr_K} under the least favorable model which is solution of the mini-max problem (\ref{minimax_pb}), i.e. the centralized problem. The performance assessment is given by the mean and variance of the least favorable state prediction error for each node $k$ (including the diffusion step), say $\tilde x_{k,t}$ with $k=1\ldots N$.
In \cite{ROBUST_STATE_SPACE_LEVY_NIKOUKHAH_2013,STATETAU_2017} it has been shown that the least favorable model can be characterized over a finte interval $[0,T]$ and it takes the following form:
\alg{\label{fl_model}\xi_{t+1}&=\check A_t\xi_t+\check B_t\varepsilon_t\nn\\
y_t&=\check C_t \xi_t+\check D_t\varepsilon_t}
where $\xi_t=[\,x_t^T\; e_t^T\,]^T$, $x_t$ is the least favorable state process, $e_t$ is the least favorable prediction error of $x_t$ using the robust filter (\ref{rob_est}) and  $\varepsilon_t$ is WGN with covariance matrix equal to the identity. Moreover,
\alg{\label{ls_model}\check A_t &:= \left[\begin{array}{cc} A &  \Gamma_B \Gamma_{H_t} \\  0&  A-G_t C+(\Gamma_B-G_t\Gamma_D)\Gamma_{H_t} \\ \end{array}\right]\nn\\
\check B_t&:= \left[\begin{array}{cc} \Gamma_B\Gamma_{L_t}  \\    (\Gamma_B-G_t\Gamma_D)\Gamma_{L_t} \\ \end{array}\right]   \nn\\
\check C_t &:= \left[\begin{array}{cc} C &  \Gamma_D \Gamma_{H_t} \\ \end{array}\right],\;\;
\check D_t := \Gamma_D \Gamma_{L_t}, }
where
$\Gamma_{L_t}$ is such that $ K_t=\Gamma_{L_t}\Gamma_{L_t}^T$,
\alg{ K_{t}&:=(I-(\Gamma_B-G_t\Gamma_D)^T (\Omega_{t+1}^{-1}+\theta_t I) (\Gamma_B-G_t\Gamma_D))^{-1} \nn\\
\Gamma_{H_t}&:= K_{t} (\Gamma_B-G_t\Gamma_D)^T  (\Omega_{t+1}^{-1}+\theta_t I)(A-G_tC).\nn}
Matrix $\Omega_{t+1}^{-1}$ is computed from the backward recursion
\alg{\label{recursion_in_omega}\Omega_t^{-1}=(A & -G_tC)^T(\Omega_{t+1}^{-1}+\theta_t I) (A-G_tC)\nn\\ &  +\Gamma_{H_t} ^T  K_{t}^{-1} \Gamma_{H_t}}
where the final point is initialized with $\Omega_{T+1}^{-1}=0$ and $T$ is the simulation horizon. Therefore, to construct the least favorable model we need to compute the gain $G_t$ performing a forward sweep of the robust Kalman filter in (\ref{rob_est}) over the interval $[0,T]$, then we generate the matrices $\Omega_t$ through a backward sweep over the interval $[0,T]$. We partition $\Gamma_{H_t}\in\Rs^{(pN+n) \times n}$ and $\Gamma_{L_t}\in\Rs^{(pN +n) \times (pN+n)}$ as $\Gamma_{H_t}=[\, M_t^T \; H_t^T\, ]^T$ and $\Gamma_{L_t}=[\, N_t^T \; L_t^T\, ]^T$ with $M_t\in\Rs^{n \times n}$, $H_t\in\Rs^{pN\times n}$, $N_t\in\Rs^{n \times (pN+n)}$ and $L_t\in\Rs^{pN \times (pN+n)}$. Moreover, we partition $H_t$ and $L_t$ as follows:   
\alg{H_t=\left[\begin{array}{c} H_{1,t} \\ \vdots \\ H_{N,t}\end{array}\right], \; \; L_t=\left[\begin{array}{c} L_{1,t} \\ \vdots \\ L_{N,t}\end{array}\right] }
where $H_{k,t}\in\Rs^{p\times n}$ and $L_{k,t}\in\Rs^{p\times (pN+n)}$.

Next we express the least favorable state prediction error $\tilde x_{k,t}$ at node $k$ in terms of the WGN $\varepsilon_t$: in this way we will be able to characterize the mean and the variance of $\tilde x_{k,t}$.  We define 
\alg{\tilde x_{k,t}=x_t-\hat x_{t,k}\nn\\
\tilde \psi_{k,t}=x_t-\psi_{k,t}}
which represent the prediction error and the intermediate prediction error, respectively, at node $k$ at time $t$. Notice that 
\alg{y_{l,t}= [\,C_l\;\, D_l H_{l,t}\,]\xi_t+D_lL_{l,t}\varepsilon_t} for $l=1\ldots N$.

 Therefore, we have
\alg{	\label{form1}\tilde \psi_{k,t+1}&=x_{t+1}-\psi_{k,t+1}=[\, I \;\, 0\,]\xi_{t+1}-\psi_{k,t+1}\nn\\
&=[\, I \;\, 0\,](\check A_t\xi_t+\check B_t \varepsilon_t)-A\hat x_{k,t}\nn\\ &\hspace{0.5cm}-A(V_{k,t}^{-1}+S_k)^{-1}\sum_{l\in\Nc_k}
C_l^T R_l^{-1} (y_{l,t}-C_l\hat x_{k,t})\nn\\
&=  A x_t+\Gamma_B\Gamma_{H_t}e_t+\Gamma_B \Gamma_{L_t}\varepsilon_t
-A\hat x_{k,t}\nn\\ &\hspace{0.5cm}-A(V_{k,t}^{-1}+S_k)^{-1}\sum_{l\in\Nc_k}
C_l^T R_l^{-1} ([\,C_l\;\, D_l H_{l,t}\,]\xi_t\nn\\ & \hspace{0.5cm}+D_lL_{l,t}\varepsilon_t-C_l\hat x_{k,t})\nn\\
&=  A \tilde x_{k,t}+B M_te_t+B N_t\varepsilon_t
\nn\\ &\hspace{0.5cm}-A(V_{k,t}^{-1}+S_k)^{-1}\sum_{l\in\Nc_k}
C_l^T R_l^{-1} (C_l x_t+D_l H_{l,t}e_t\nn\\ & \hspace{0.5cm}+D_lL_{l,t}\varepsilon_t-C_l\hat x_{k,t})\nn\\
&=  A \tilde x_{k,t}+BM_te_t+B N_t\varepsilon_t
\nn\\ &\hspace{0.5cm}-A(V_{k,t}^{-1}+S_k)^{-1}\sum_{l\in\Nc_k}
C_l^T R_l^{-1} (C_l \tilde x_{k,t}+D_l H_{l,t}e_t\nn\\ & \hspace{0.5cm}+D_lL_{l,t}\varepsilon_t)\nn\\
&=  A(I- (V_{k,t}^{-1}+S_k)^{-1}S_k) \tilde x_{k,t}+B(M_te_t+N_t\varepsilon_t)
\nn\\ &\hspace{0.5cm}-A(V_{k,t}^{-1}+S_k)^{-1}\sum_{l\in\Nc_k}
C_l^T R_l^{-1} D_l (H_{l,t}e_t +L_{l,t}\varepsilon_t)\nn\\
&=  A(V_{k,t}^{-1}+S_k)^{-1}V_{k,t}^{-1} \tilde x_{k,t}+B(M_te_t+N_t\varepsilon_t)
\nn\\ &\hspace{0.5cm}-A(V_{k,t}^{-1}+S_k)^{-1}\sum_{l\in\Nc_k}
C_l^T R_l^{-1} D_l (H_{l,t}e_t +L_{l,t}\varepsilon_t)\nn\\
&=  A(V_{k,t}^{-1}+S_k)^{-1}V_{k,t}^{-1} \tilde x_{k,t}+B(M_te_t+N_t\varepsilon_t)
\nn\\ &\hspace{0.25cm}-A(V_{k,t}^{-1}+S_k)^{-1}\sum_{l=1}^N j_{lk}
C_l^T R_l^{-1} D_l (H_{l,t}e_t +L_{l,t}\varepsilon_t) }
where we recall that $j_{lk}=1$ if $l\in\Nc_k$, otherwise $j_{lk}=0$. Notice that,
\alg{\label{form2}\tilde x_{k,t+1}&=x_t-\sum_{l\in\Nc_k} w_{lk} \psi_{l,t+1}\nn\\ &= \sum_{l=1}^N w_{lk} x_t-\sum_{l=1}^N w_{lk} \psi_{l,t+1}\nn\\&=\sum_{l=1}^N w_{lk} (x_t-\psi_{l,t+1})=\sum_{l=1}^N w_{lk} \tilde \psi_{l,t+1}}
where we have exploited (\ref{prop_W}).
Taking into account (\ref{form1}), we have 
\alg{\label{form3} \tilde x_{k,t+1}&=\sum_{l=1}^N w_{lk} [ 
  A(V_{l,t}^{-1}+S_l)^{-1}V_{l,t}^{-1} \tilde x_{l,t} -A(V_{l,t}^{-1}+S_l)^{-1}\nn\\  &\hspace{0.5cm}\times \sum_{m=1}^N j_{ml}
C_m^T R_m^{-1} D_m (H_{m,t}e_t +L_{m,t}\varepsilon_t)]\nn\\ &\hspace{0.5cm}+B(M_te_t+N_t\varepsilon_t).} By defining
\alg{\tilde \chi_t &:= [\, \tilde x_{1,t}^T \; \ldots \; \tilde x_{N,t}^T\,]^T\nn\\
\Cc&:=\mathrm{diag}(C_{1},\ldots ,C_{N})\nn\\
 \Vc_t&:=\mathrm{diag}(V_{1,t},\ldots ,V_{N,t})\nn\\
\Sc&:=\mathrm{diag}(S_{1},\ldots ,S_{N})
} 
we can rewrite (\ref{form3}) in the following compact way:
\alg{\tilde x_{k,t+1}&=([\,w_{1 k}\; \ldots \; w_{Nk} \,]\otimes I)\{(I\otimes A)(\Vc_t^{-1}+\Sc)^{-1}\nn\\ &  \hspace{0.5cm}\times\Vc^{-1}_t \tilde \chi_t
-(I\otimes A)(\Vc_t^{-1}+\Sc)^{-1}(J^T\otimes I)\nn\\ 
& \hspace{0.5cm} \times \Cc^T R^{-1} D(H_te_t+L_t\varepsilon_t)\}+B(M_te_t+N_t\varepsilon_t)\nn}
and thus
\alg{\tilde \chi_{t+1}&=(W^T\otimes I)\{(I\otimes A)(\Vc_t^{-1}+\Sc)^{-1}\Vc^{-1}_t \tilde \chi_t\nn\\ &  \hspace{0.5cm}
-(I\otimes A)(\Vc_t^{-1}+\Sc)^{-1}(J^T\otimes I )\Cc^T R^{-1} D\nn\\ 
& \hspace{0.5cm} \times (H_te_t+L_t\varepsilon_t)\}+ (\mathbf{1} \otimes B)({M_t}e_t+{N_t}\varepsilon_t)}
where $\mathbf 1$ denotes the vector of ones. We rewrite the latter as 
\alg{\label{form4}\tilde \chi_{t+1} &= \Ac_t \tilde \chi_t+\Bc_t \varepsilon_t+\Cc_t e_t }
where 
\alg{\Ac_t&:= (W^T \otimes I)(I\otimes A)(\Vc^{-1}_t+\Sc)^{-1}\Vc_t^{-1}\nn\\
\Bc_t&:= -(W^T\otimes I)(I\otimes A)(\Vc_t^{-1}+\Sc)^{-1}(J^T \otimes I )\Cc R^{-1} D L_t \nn\\ &\hspace{0.5cm}+\mathbf{1}\otimes B N_t\nn\\
\Cc_t&:= -(W^T\otimes I)(I\otimes A)(\Vc_t^{-1}+\Sc)^{-1}(J^T \otimes I )\Cc R^{-1} D H_t \nn\\ &\hspace{0.5cm}+\mathbf{1}\otimes B M_t.}
Combining (\ref{form4}) with the model for $e_t$ in (\ref{fl_model}), we obtain:
\alg{\label{form5}\eta_{t+1}=\Fc_t\eta_t+\Gc_t\varepsilon_t}
where $\eta_t:=[\, \tilde \chi_t^T\; e_t^T\,]^T$,
\alg{\Fc_t&:=\left[\begin{array}{cc}\Ac_t  & \Cc_t \\ 0  & (A-G_tC)+ (\Gamma_{B}-G_t\Gamma_{D}) \Gamma_{H_t}\\ \end{array}\right]\nn\\  \Gc_t&:=\left[\begin{array}{c}\Bc_t \\ (\Gamma_{B}-G_t\Gamma_{D}) \Gamma_{L_t}\\ \end{array}\right].}

From (\ref{form5}), we can analyze the performance of the distributed Algorithm \ref{rob_est_diff}. Taking the expectation of (\ref{form5}), we obtain 
\alg{\tilde \Es[\eta_{t+1}]=\Fc_t \tilde \Es[\eta_t].}
Since $\hat x_0$ is the mean of $x_0$ and $\hat x_{k,0}=\hat x_0$ for $k=1\ldots N$, we have that $\tilde \Es[\eta_0]=0$. Accordingly, $\eta_t$ is a random vector with zero mean for any $t$. This means that the distributed Kalman predictions   
with diffusion step are unbiased.

We proceed to analyze the variance of the prediction errors. 
We define $\Qc_t=\tilde \Es[\eta_t\eta_t^T]$. Since $\varepsilon_t$ is WGN with covariance matrix equal to the identity, by (\ref{form5}) we have that $\Qc_t$ is given by solving the following Lyapunov equation
\alg{\label{Layp_eq}\Qc_{t+1}=\Fc_t\Qc_t\Fc_t^T+\Gc_t\Gc_t^T.}  
We partition $\Qc_t$ as follows:
\alg{\Qc_t=\left[\begin{array}{ccc}\Pc_t  & \Hc_t \\ \Hc_t^T & \Rc_t \end{array}\right]}
where $\Pc_t\in\Rs^{Np\times Np}$, $\Hc_t\in\Rs^{Np\times n}$ and $\Rc_t\in\Rs^{n\times n}$. The $n\times n$ matrices in the  main block diagonal of $\Pc_t$ represent the covariance matrices of the estimation error at each node. 
Let $\mathrm{MSD}_{k,t}:=\tilde \Es[\|x_t-\hat x_{k,t}\|^2]$ denote the least favorable mean square deviation at node $k$ and at time $t$. Then, the average least favorable mean square deviation across the network at time $t$ is 
\alg{\overline{\mathrm{MSD}}_t:=\frac{1}{N}\sum_{k=1}^N \mathrm{MSD}_{k,t}=\frac{1}{N}\tr(\Pc_t).}
The computation of the sequence $\Pc_t$ depends on the simulation horizon $T$. In particular, it is required to perform three steps:
\begin{itemize}
\item compute the filtering gain $G_t$ performing a forward sweep of the centralized robust Kalman filter in (\ref{rob_est}) over the interval $[0,T]$
\item compute $\Omega_{t}$ performing the backward recursion (\ref{recursion_in_omega}) over the interval $[0,T]$
\item compute $\Pc_t$ performing a forward sweep of the Lyapunov equation in (\ref{Layp_eq}) over the interval $[0,T]$.
\end{itemize}

\subsection{Convergence analysis}
In the previous section we showed how to compute $\Qc_t$ over the simulation horizon $[0,T]$. Let $0<\alpha<\beta<1$. We show that under reachability and local observability, and choosing the tolerance $c>0$ sufficiently small, then $\Qc_t$ converges over the interval $[\alpha T,\beta T]$ as $T$ approaches infinity, and thus the prediction errors at each node have zero mean and finite constant variance in steady state. It is worth noting that local observability is a strong assumption which may pose some limits on the practical applicability of the distributed robust Kalman filter.

In \cite{zorzi2018lf} it has been shown that the least favorable prediction error using an estimator of type $\hat x_{t+1}=A\hat x_t+\tilde G_t(y_t-C\hat x_t)$, where $C$ coincides with the one in the least favorable model,  has zero mean and convergent covariance matrix provided that $c$ is sufficiently small. However, these results cannot be directly applied to our case because the predictor at node $k$ is given by a convex combination of local estimators whose matrix $C$ does not coincide with the one of the least favorable model.

The update of the intermediate local prediction can be rewritten as
\alg{\psi_{k,t+1}=A\hat x_{k,t} +G_{k,t}(y_{k,t}^{loc}-C_{k}^{loc} \hat x_{k,t}) }
where 
\alg{G_{k,t}=A(V_{k,t}^{-1}+S_k^{-1})^{-1} (C_k^{loc})^T(R_k^{loc})^{-1}} is the filtering gain at node $k$. The first step is to show that $G_{k,t}$ converges as $t$ approaches infinity.

\begin{propo}\label{prop_conv_filtro}Assume that the pair $(A,B)$ is reachable and that the pair $(A,C_k^{loc})$ is observable for every $k$. Then, there exists $c>0$ sufficiently small such that for any $V_{k,0}>0$ the sequence $P_{k,t}$ $t\geq 0$ generated by Algorithm \ref{rob_est_diff} at node $k$ converges to a unique solution $\bar P_k>0$. Furthermore, $\theta_{k,t}\rightarrow \bar \theta_k$, $V_{k,t}\rightarrow \bar V_k>0$ and the limit $\bar G_k$ of the filtering gain $G_{k,t}$  is such that $A-\bar G_kC_k^{loc}$ is Schur stable. Moreover, $\bar P_k$ is the unique solution of the algebraic Riccati-like equation
\alg{\bar P_k=A&( \bar P_k^{-1}-	\bar\theta_k I+(C_k^{loc})^T(R_k^{loc})^{-1}C_k^{loc})^{-1}A^T+B B^T.}
  \end{propo}

\IEEEproof The convergence of the local robust Kalman filter follows from the convergence result of the robust Kalman filter in \cite[Proposition 3.5]{ZORZI_CONTRACTION_CDC}, see also \cite{CONVTAU}, under the assumption that the local state space model $(A,B,C_k^{loc},D_k^{loc})$ is reachable and observable. \qed

Regarding the least favorable model in (\ref{fl_model}), it is possible to prove that it does converge to a state space model with constant parameters, as the simulation horizon $T$ tends to infinity.
\begin{propo}[Zorzi, Levy \cite{zorzi2018lf}] \label{propp_lf_CDC}Assume that the pair $(A,B)$ is reachable and that the pair $(A,C)$ is observable. Then, there exists $c>0$ sufficiently small such that:
\begin{itemize}
\item the forward sequences $G_t$ and $\theta_t$, $t\geq 0$, of the centralized robust Kalman filter (\ref{rob_est}) converges to $\bar G$ and $\bar \theta$, respectively, as $t$ tends to infinity;
\item when the simulation horizon $T$ tends to infinity, the backward sequence $\Omega_t$ generated by (\ref{recursion_in_omega}), with the steady state parameters $\bar G$ and $\bar \theta$ of the centralized robust Kalman filter, converges to $\bar \Omega$. Furthermore,  
\alg{L_t \rightarrow \bar L, \hspace{0.7cm}  K_t \rightarrow \bar K,  \hspace{0.7cm}
\Gamma_{H_t}\rightarrow  \Gamma_{\bar H}
} and $(A-\bar GC)+(B- \bar G\Gamma_{D})\Gamma_{\bar H}$ is Schur stable.
\end{itemize}
\end{propo}

Finally, we need of the following result. \begin{propo}[Cattivelli, Sayed \cite{cattivelli2010}] \label{prop_Lyap_catt} Consider the time-varying Lyapunov equation \alg{X_{t+1}=A_tX_tA_t+Q_t}
where $A_t$ and $Q_t$ converges to $A$ and $Q$, respectively, as $t\rightarrow \infty$, with $A$ Schur stable. Then, $X_t$ converges to the unique solution $X$ of the Lyapunov equation:
\alg{X=AXA^T+Q.}
\end{propo}

We are ready to prove the main convergence result. 
\begin{propo}Assume that the pair $(A,B)$ is reachable and that the pair $(A,C_k^{loc})$ is observable for every $k$. Then, there exists $c>0$ sufficiently small such that for any $V_0>0$ and $V_{k,0}>0$ the sequence $\Qc_{t}$, $t\geq 0$, generated by (\ref{Layp_eq}) converges to $\bar \Qc>0$, $\Fc_t\rightarrow \bar \Fc$ and $\Gc_t\rightarrow \bar\Gc$ 
over $[\alpha T, \beta T]$ as $T\rightarrow \infty$. Moreover, $\bar \Qc$ is the unique solution of
\alg{\label{Layp_eq_alg}\bar \Qc =\bar \Fc \bar \Qc \bar\Fc ^T+\bar \Gc \bar \Gc ^T} and $\bar\Fc$ is Schur stable. Therefore, the average least favorable mean square deviation across the network $\overline{\mathrm{MSD}}_t$ does converge over $[\alpha T, \beta T]$ as $T\rightarrow \infty$.
\end{propo}
\IEEEproof
Notice that the assumptions of Proposition \ref{prop_conv_filtro} hold. Therefore, we have that $\Vc_t\rightarrow \bar \Vc $ and $\Ac_t\rightarrow \bar \Ac$ where \alg{\bar\Vc&:=\mathrm{diag}(\,\bar V_1,\, \ldots,\, \bar V_N \,)\nn\\
\bar \Ac&:=(W^T\otimes I)(I\otimes A)(\bar \Vc^{-1}+\Sc)^{-1}\bar\Vc^{-1}.} Moreover, $\bar G_{k}=A(\bar V_{k}^{-1}+S_k^{-1})^{-1} (C_k^{loc})^T(R_k^{loc})^{-1}$ and thus 
\alg{\label{formAG_k} A- &\bar G_k C_k^{loc}\nn\\ &=A(I-(\bar V_k^{-1}+  S_k)^{-1}(C_k^{loc})^T(R_k^{loc})^{-1} C_k^{loc})\nn\\ &= A(I-(\bar V_k^{-1}+  S_k)^{-1}S_k)\nn\\ &=A(\bar V_k^{-1} +S_k)^{-1} \bar V_k^{-1}} which is Schur stable. Accordingly, we have that the block-diagonal 
matrix \alg{\mathcal M:=(I\otimes A)(\bar \Vc^{-1}+\Sc)^{-1}\bar\Vc^{-1} } is Schur stable. Then, by using \cite[Lemma 2]{cattivelli2010} we have that $\bar \Ac=(W^T\otimes I)\mathcal M$ is Schur stable, because $W$ satisfies the conditions in (\ref{prop_W}).

Since the assumptions of Proposition \ref{propp_lf_CDC} hold, then $\Fc_t\rightarrow \bar \Fc$ and $\Gc_t\rightarrow \bar \Gc$ over the interval $[\alpha T, \beta T]$ as $T\rightarrow \infty$. Moreover,
\alg{\bar \Fc&=\left[\begin{array}{cc}  \bar \Ac & \star \\0   & (A-\bar GC)+(\Gamma_B-\bar G\Gamma_D)\Gamma_{\bar H}\end{array}\right]\nn\\
\bar \Gc&=\left[\begin{array}{c}  \bar \Bc  \\ (\Gamma_B-\bar G\Gamma_D)\Gamma_{\bar L}\end{array}\right]\nn} 
where $\bar\Bc:= -(W^T\otimes I)(I\otimes A)(\bar \Vc^{-1}+\Sc)^{-1}(J^T \otimes I )\Cc R^{-1} D \bar L +\mathbf{1}\otimes B \bar N$
with $\Gamma_{L_t}\rightarrow \Gamma_{\bar L}=[\,\bar N^T\; \bar L^T\,]^T$. Since $\bar \Fc$ is block upper-triangular, its eigenvalues coincides with the eigenvalues of  $\bar \Ac$ and $(A-\bar GC)+(\Gamma_B-\bar G\Gamma_D)\Gamma_{\bar H}$. Since the latter are Schur stable, we conclude that $\bar \Fc$ is Schur stable.

The Lyapunov equation in (\ref{Layp_eq}) satisfies the assumptions of Proposition \ref{prop_Lyap_catt} for $c$ sufficiently small. Accordingly, $\Pc_t$ converges to $\bar \Pc$ over $[\alpha T, \beta T]$ as $T\rightarrow \infty$ which is the unique solution to (\ref{Layp_eq_alg}). \qed 
 
\subsection{Optimality property under large deviations} \label{sec_opt}
We shall show that, under the least favorable model solution to (\ref{minimax_pb}) with $c$ sufficiently large, it is very likely that the predictor at node $k$ of Algorithm \ref{rob_est_diff} performs better than the predictor at node $k$ based on the scheme in \cite{cattivelli2010}. The latter scheme, indeed, does not consider the possibility that the actual model and  the nominal model do not coincide exactly.

Without loss of generality we assume that $y_t$ is obtained by stacking first $y_{l,t}$ with $l\in \Nc_k$ and then $y_{l,t}$ with $l\notin \Nc_k$. Then, it is not difficult to see that that the standard local Kalman predictor at node $k$ coincides with the Kalman predictor based on the model
\alg{\label{modello_phi_loc}  x_{t+1}&=Ax_t+B w_t\nn\\
y_t&= \left[\begin{array}{c} C^{loc}_k\\ 0\end{array}\right]x_t+ \left[\begin{array}{c} D^{loc}_k\\ \breve F_{k,t}^{loc}\end{array}\right]v_t}
where $\breve F_{k,t}^{loc}$ is an arbitrary invertible matrix because the observations $y_{l,t}$ with $l\notin \Nc_k$ play no role at node $k$. Let  $Q_{k,t}^{loc}:=\breve F_{k,t}^{loc}(\breve F_{k.t}^{loc})^T$. Therefore the robust intermediate prediction at node $k$ in Algorithm \ref{rob_est_diff}, hereafter denoted by \textsf{RKF diff}, is the solution of the mini-max problem
\alg{\label{minimax_loc}\hat \psi_{k,t+1} &= \underset{g_t\in\mathbb G_{k,t}}{\mathrm{argmin} \,}\underset{\tilde \phi_{k,t}\in\mathbb B_{k,t}}{\mathrm{max}} \tilde \Es_k[\| x_{t+1}-g_t\|^2|Y_{t-1}]}
where
\alg{\mathbb{B}_{k,t}:=\{\, \tilde \phi_{k,t} \hbox{ s.t. } \tilde\Es_k[\log (\tilde \phi_{k,t} / \phi_{k,t})| Y_{t-1}]\leq c\,\};}
$\phi_{k,t}$ is the transition probability density of $z_t$ given $x_t$ corresponding to (\ref{modello_phi_loc}) and $\tilde \phi_{k,t}$ is the least favorable one in $\mathbb B_{k,t}$; $\tilde\Es_k[\log (\tilde \phi_{k,t} / \phi_{k,t})| Y_{t-1}]$ and $\tilde \Es_k[\| x_{t+1}-g_t\|^2|Y_{t-1}]$
are defined as in (\ref{def_ball_E}) and (\ref{minimax_pbE}), respectively, with $\phi_t$, $\tilde  \phi_t$ and $\tilde f_t$ replaced by 
$ \phi_{k,t}$, $\tilde \phi_{k,t}$ and $\tilde f_{k,t}$, respectively; $\tilde f_{k,t}(x_t|Y_{t-1})\sim \Nc(\hat x_{k,t},V_{k,t})$ is the least favorable conditional probability density of $x_t$ given $Y_{t-1}$ at node $k$; $\mathbb G_{k,t}$ denotes the set of all estimators such that $\tilde\Es_k[\|g_t\|^2]$ is finite for any $\tilde \phi_{k,t}\in\mathbb B_{k,t}$.

In 
\cite{ROBUST_STATE_SPACE_LEVY_NIKOUKHAH_2013} it has been shown that the mini-max problem (\ref{minimax_loc}) can be reformulated in terms of
\alg{ \bar p_t^{loc}(z_t|Y_{t-1})&=\int \phi_{k,t}(z_t|x_{t})\tilde f_{k,t}(x_t|Y_{t-1})\mathrm d x_t\\
\tilde p_t^{loc}(z_t|Y_{t-1})&=\int \tilde \phi_{k,t}(z_t|x_{t})\tilde f_{k,t}(x_t|Y_{t-1})\mathrm d x_t}
representing the local pseudo-nominal and the local least favorable conditional probability density of $z_t$ given $Y_{t-1}$. Moreover, $\mathbb D_{KL}(\tilde p_t^{loc}, \bar p_t^{loc})=c$ and
\alg{\mathbb D_{KL}(\tilde p_t^{loc},\bar p_t^{loc})=\int \tilde p_t^{loc}(z_t|Y_{t-1})\log\left(\frac{\tilde p_t^{loc}(z_t|Y_{t-1})}{\bar p_t^{loc}(z_t|Y_{t-1})}\right)\mathrm d z_t}
is the Kullback-Leibler divergence between $\tilde p_t^{loc}$ and $\bar p_t^{loc}$.
\begin{rem} \label{rem_lik}It is worth noting that $\mathbb D_{KL}( p_1,p_2)$ represent the negative log-likelihood (up to constant factors) of $p_2$ under the the model described by $  p_1$, \cite{robert2014machine}. Assume that $\mathbb D_{KL}( p_1,p_2)\ll \mathbb D_{KL}( p_1,p_3)$ where the symbol $\ll$ means ``much less than''. This means that $p_2$ explains the data generated by $p_1$ better than $p_3$.\end{rem}

Notice that the best intermediate prediction at node $k$ is the one constructed using the least favorable model 
\alg{ \tilde p_t(z_t|Y_{t-1})&=\int \tilde \phi_{t}(z_t|x_{t})\tilde f_{k,t}(x_t|Y_{t-1})\mathrm d x_t }
since we are are evaluating the performance under the least favorable model solution to the ``centralized'' mini-max problem in (\ref{minimax_pb}) which is not available because it requires to compute the centralized filtering gains of the centralized robust Kalman filter. On the other hand, 
the intermediate prediction at node $k$ of the algorithm proposed in \cite{cattivelli2010}, hereafter denoted by \textsf{KF diff}, is constructed using the nominal model
\alg{ p_t^{loc}(z_t|Y_{t-1})&=\int  \phi_{k,t}(z_t|x_{t}) f_{k,t}(x_t|Y_{t-1})\mathrm d x_t} where $f_{k,t}(x_t|Y_{t-1})$ is the nominal conditional probability density of $x_t$ given $Y_{t-1}$ at node $k$.
 In view of Remark \ref{rem_lik}, the next proposition shows that if $c$ is sufficiently large, then $\tilde p_t^{loc}$ explains the measurements generated by the actual model $\tilde p_t$ better than $p_t^{loc}$. As a consequence, it is very likely that the performance of $\psi_{k,t+1}$ using  \textsf{RKF diff} (i.e. using $\tilde p_t^{loc}$) is better than the one using  \textsf{KF diff} (i.e. using $p_t^{loc}$).

\begin{propo}Assume that for some $t$ the distribution of $x_{k,t}$ given $Y_{t-1}$ at node $k$ is fixed and it is the same for \textsf{RKF diff}, \textsf{KF diff}, that is $f_{k,t}(x_t|Y_{t-1})$ and $\tilde f_{k,t}(x_t|Y_{t-1})$ coincide. Then, for $c$ sufficiently large we have that 
\alg{\label{claim_opt} \mathbb D_{KL}(\tilde p_t,\tilde p_t^{loc} ) \ll \mathbb D_{KL} (\tilde p_t,p_t^{loc}).}
\end{propo}
\IEEEproof
Let $\tilde f_{k,t}\sim \Nc(\hat x_{k,t},V_{k,t})$ with $V_{k,t}>0$ which is fixed and thus it does not depend on $c$. First, notice that $p_t(z_t|Y_{t-1})=\bar p_t^{loc}(z_t|Y_{t-1})$ because the distribution of $x_{k,t}$ given $Y_{t-1}$ is the same for \textsf{RKF diff} and \textsf{KF diff}.
Accordingly,
\alg{ \label{pdfs_gauss}p_t^{loc}(z_t|Y_{t-1}) &\sim \Nc(\mu_t^{loc},K_t^{loc})\nn \\
\tilde p_t^{loc}(z_t|Y_{t-1}) &\sim \Nc(\mu_t^{loc},\tilde K_t^{loc})\nn \\
\tilde p_t(z_t|Y_{t-1}) &\sim \Nc(\mu_t,\tilde K_t)}
where 
\alg{\mu_t^{loc}=\left[\begin{array}{c} A\\ C_k^{loc}\\ 0\end{array}\right] \hat x_{k,t},\; \; \mu_t=\left[\begin{array}{c} A\\ C_k^{loc}\\ \breve C_k^{loc}\end{array}\right]\hat x_{k,t},}
{\small  \alg{  K_t^{loc}&=\left[\hspace{-0.2cm}\begin{array}{c} A\\ C_k^{loc} \\ 0\end{array}\hspace{-0.2cm}\right]V_{k,t} \left[\hspace{-0.2cm} \begin{array}{ccc}A^T &  \hspace{-0.2cm}(C_k^{loc})^T &\hspace{-0.2cm} 0  \end{array}\hspace{-0.2cm} \right]+\left[\hspace{-0.2cm}\begin{array}{ccc}BB^T &  \hspace{-0.2cm}0&  \hspace{-0.2cm} 0\\ 0 & \hspace{-0.2cm} R_k^{loc} & \hspace{-0.2cm} 0\\ 0 &  \hspace{-0.2cm} 0 & \hspace{-0.2cm} Q_{k,t}^{loc}  \end{array}\hspace{-0.2cm}\right],\nn\\
\tilde  K_t^{loc}&=  K_t^{loc}+\left[\hspace{-0.2cm}\begin{array}{c} I\\ 0 \\ 0\end{array}\hspace{-0.2cm}\right](V_{k,t+1} -P_{k,t+1} ) \left[\hspace{-0.2cm} \begin{array}{ccc}I &  \hspace{-0.2cm} 0&\hspace{-0.2cm} 0  \end{array}\hspace{-0.2cm} \right],\nn\\
 K_t&=\left[\hspace{-0.2cm}\begin{array}{c} A\\ C_k^{loc} \\ \breve C_k^{loc} \end{array}\hspace{-0.2cm}\right]V_{k,t} \left[\hspace{-0.2cm} \begin{array}{ccc}A^T &  \hspace{-0.2cm}(C_k^{loc})^T &\hspace{-0.2cm} (\breve C_k^{loc})^T  \end{array}\hspace{-0.2cm} \right]+\left[\hspace{-0.2cm}\begin{array}{ccc}BB^T & \hspace{-0.2cm}  0& \hspace{-0.2cm} 0\\ 0 & \hspace{-0.2cm} R_k^{loc} &  \hspace{-0.2cm} 0\\ 0 &  \hspace{-0.2cm} 0 & \hspace{-0.2cm} \breve R_k^{loc}  \end{array}\hspace{-0.2cm}\right],\nn\\
 \tilde K_t&=K_t+\left[\hspace{-0.2cm}\begin{array}{c} I\\ 0 \\ 0\end{array}\hspace{-0.2cm}\right](V_{t+1}-P_{t+1}) \left[\hspace{-0.2cm} \begin{array}{ccc} I&  \hspace{-0.2cm}0 &\hspace{-0.2cm}  0  \end{array}\hspace{-0.2cm} \right];\nn
 } }
$\breve C_k^{loc}$ and $\breve R_k^{loc}$ are the matrices obtained by using $C_l$ and $R_l$, respectively, with $l\notin\Nc_k $. It is worth noting that the relation between $\tilde K_t^{loc}$ and $K_t^{loc}$ is given by \cite[Theorem 1]{ROBUST_STATE_SPACE_LEVY_NIKOUKHAH_2013}. The same observation holds between $\tilde K_t$ and $K_t$ where the latter represents the covariance matrix of $z_t$ given $Y_{t-1}$ in the nominal model. Moreover,
{\small \alg{P_{k,t+1}&=A V_{k,t}A^T-A V_{k,t}  (C_k^{loc})^T  \left(  C_k^{loc}   V_{k,t}   (C_k^{loc})^T  + R_k^{loc}  \right)^{-1}  C_k^{loc}  V_{k,t}A^T+BB^T\nn\\
V_{k,t+1}&=(P_{k,t+1}^{-1}-\theta_{k,t}I)^{-1}\nn\\
P_{t+1}&=A V_{k,t}A^T-A V_{k,t}  \left[\hspace{-0.1cm} \begin{array}{cc}   \hspace{-0.1cm}(C_k^{loc})^T &\hspace{-0.1cm} (\breve C_k^{loc})^T  \end{array}\hspace{-0.2cm} \right] \left( \left[\hspace{-0.2cm}\begin{array}{c} C_k^{loc} \\ \breve C_k^{loc} \end{array}\hspace{-0.2cm}\right]  V_{k,t}  \left[\hspace{-0.1cm} \begin{array}{cc}   \hspace{-0.1cm}(C_k^{loc})^T &\hspace{-0.1cm} (\breve C_k^{loc})^T  \end{array}\hspace{-0.2cm} \right]+\left[\hspace{-0.1cm}\begin{array}{cc}  \hspace{-0.1cm} R_k^{loc} &  \hspace{-0.1cm} 0\\    \hspace{-0.1cm} 0 & \hspace{-0.1cm} \breve R_k^{loc}  \end{array}\hspace{-0.2cm}\right] \right)^{-1} \left[\hspace{-0.2cm}\begin{array}{c} C_k^{loc} \\ \breve C_k^{loc} \end{array}\hspace{-0.2cm}\right]  V_{k,t}A^T+BB^T\nn\\
V_{t+1}&=(P_{t+1}^{-1}-\theta_{t}I)^{-1}\nn 
} }
and $\theta_{k,t}$, $\theta_t$ are the solution to $\gamma(P_{k,t+1},\theta_{k,t})=c$,  $\gamma(P_{t+1},\theta_{t})=c$, respectively. Recall that 
 \alg{\label{def_gamma_bis}\gamma(P,\theta):=\log\det(I-\theta P)+\tr((I-\theta P)^{-1}-I).}

In view of (\ref{pdfs_gauss}), it is not difficult to see that 
\alg{\mathbb D_{KL}(\tilde p_t, \tilde p_t^{loc})=\mathbb D_{KL}(\tilde p_t,  p_t^{loc})+\frac{1}{2}d_\Delta} where \alg{d_\Delta &= \delta^T ((\tilde K_t^{loc})^{-1}-( K_t^{loc})^{-1} )\delta+\log \det(\tilde K_t^{loc})  \nn\\
& \hspace{0.35cm}-\tr( \tilde K_t (K_t^{loc})^{-1} )+\tr(\tilde K_t(\tilde K_t^{loc})^{-1})  -\log \det( K_t^{loc}) \nn\\
& \leq \log \det(\tilde K_t^{loc}) +\tr\left[\tilde K_t \left((\tilde K_t^{loc})^{-1}-(K_t^{loc})^{-1}\right) \right]   -\log \det( K_t^{loc})}
where $\delta=\mu_t-\mu_t^{loc}$ and we have exploited the fact that $(\tilde K_t^{loc})^{-1}-( K_t^{loc})^{-1}	\leq 0$ because $P_{k,t+1}<V_{k,t+1}$ and thus $\tilde K_t^{loc}\geq  K_t^{loc}$. Moreover, after some algebraic manipulations we obtain
\alg{d_\Delta\leq n\log \|V_{k,t+1}\| -\beta_{k,t} \|V_{t+1}\| +\nu_{k,t}}
where 
{\footnotesize \alg{ \beta_{k,t}&= \lambda_{min}(P_{k,t+1}^{-1} [P_{k,t+1}^{-1}+(V_{k,t+1}-P_{k,t+1})^{-1}]^{-1}P_{k,t+1}^{-1})^{-1}\tr(\bar V_{t+1}-\|V_{t+1}\|^{-1}P_{t+1} ) \nn\\
\nu_{k,t}&=-\log \det K_t^{loc}+(Np+n)\log\lambda_{max}(K_t^{loc}) + \log\det (\|V_{k,t+1}\|^{-1}I_n+\lambda_{max} (K_t^{loc})^{-1} \bar V_{k,t+1})\nn}} $\lambda_{max}(K_t^{loc})$ denotes the maximum eigenvalue of $K_t^{loc}$, $\bar V_{k,t+1}:=\|V_{k,t+1}\|^{-1}V_{k,t+1}$ and $\bar V_{t+1}:=\|V_{t+1}\|^{-1}V_{t+1}$.

It \cite{LEVY_NIKOUKHAH_2004} it has been shown that the mapping $c\mapsto \|V_{k+1,t}\|$ has singular value which is positive. Accordingly, if we take a sequence $c^{(m)}$, $m\in\mathbb N$, such that $c^{(m)}>0$ and $c^{(m)}\rightarrow \infty$ as $m\rightarrow \infty$, then $\|V_{k,t+1}^{(m)}\|\rightarrow \infty$. The same reasoning holds for the mapping $c\mapsto \|V_{t+1}\|$ and thus $\|V_{t+1}^{(m)}\|\rightarrow \infty$. Consider the sequences $\bar V^{(m)}_{k,t+1} := \| V^{(m)}_{k,t+1} \|^{-1} V_{k,t+1}^{(m)}$ and $\bar V^{(m)}_{t+1} := \| V^{(m)}_{t+1} \|^{-1} V_{t+1}^{(m)}$ which belong to the compact set $\mathcal U:=\{ \,V \hbox{ s.t. } \|V\|=1\,\}$. Therefore, there exist the subsequences $\bar V^{(m_l)}_{k,t+1}$, $l\in\mathbb N$ and $\bar V^{(m_l)}_{t+1}$, $l\in\mathbb N$, converging to $\bar V^{(\infty)}_{k,t+1}$ and $\bar V^{(\infty)}_{t+1}$, respectively. It is worth noting that $\bar V^{(\infty)}_{k,t+1},\bar V^{(\infty)}_{t+1}\geq 0$ and different from the null matrix because $\bar V^{(\infty)}_{k,t+1},\bar V^{(\infty)}_{t+1}\in\mathcal U$.  Accordingly, if we consider the corresponding subsequences for $\beta_{k,t}$ and $\nu_{k,t}$, we have: $\beta^{(m_l)}_{k,t} \rightarrow \lambda_{min}(P_{k,t+1}^{-1})^{-1}\tr(\bar V_{t+1})>0$ and $\nu_{k,t}^{(m_l)}$ is bounded above.

Next we show that $\| V_{t+1}^{(m_l)}\|/ \| V_{k,t+1}^{(m_l)}\|\rightarrow \zeta>0$. First, we recall that  $ V_{k,t+1}^{(m_l)}$ and $V_{t+1}^{(m_l)}$ are given by $\theta_{k,t}^{(m_l)}$ and $\theta_{t}^{(m_l)}$, respectively. In particular, we have $\gamma(P_{t+1}^{(m_l)},\theta_{t}^{(m_l)})=c^{(m_l)}$. Notice that we can rewrite the latter  as
\alg{\label{gamma_eig}\sum_{i=1}^n \log (1-d_i\theta_{t}^{(m_l)}) +(1-\theta d_i^{(m_l)})^{-1}-1=c^{(m_l)}}
where $d_i\geq d_{i+1}$ denotes the eigenvalues of $P_{t+1}$ and $0<\theta_{t}^{(m_l)}<d_1^{-1}$. In what follows we assume that the eigenvalue $d_1$ has multiplicity equal to one, and thus $d_1>d_i$ with $i\geq 2$. This assumption is not restrictive, indeed it generically holds. Then we can rewrite (\ref{gamma_eig}) as 
\alg{f(d_1\theta_{t}^{(m_l)})+ \breve c^{(m_l)}=c^{(m_l)}\nn}
where \alg{f(x)&=\log (1-x) +(1-x)^{-1}-1 \nn\\
\breve c^{(m_l)}&=\sum_{i=2}^n \log (1-d_i\theta_{t}^{(m_l)}) +(1-\theta d_i^{(m_l)})^{-1}-1,\nn}
$\breve c^{(m_l)} \rightarrow \breve c$ and $\breve c$ is a bounded value. Therefore
\alg{f(d_1\theta_{t}^{(m_l)})=c^{(m_l)}- \breve c^{(m_l)}.\nn} Since $c^{(m_l)}\rightarrow \infty$, we have $ \breve c^{(m_l)}=o(c^{(m_l)})$, i.e. $ \breve c^{(m_l)}/  c^{(m_l)}\rightarrow 0$ as $l$ tends to infinity. Accordingly,
\alg{\label{gamma_eig2} f(d_1\theta_{t}^{(m_l)})=c^{(m_l)}-o( c^{(m_l)}).} The same reasoning applies for $\theta_{k,t}^{(m_l)}$:
\alg{\label{gamma_eig3} f(d_{k,1}\theta_{k,t}^{(m_l)})=c^{(m_l)}-o( c^{(m_l)})}
where $d_{k,i}\geq d_{k,i+1}$ are the eigenvalues of $P_{k,t+1}$ and $d_{k,1}$ has multiplicity equal to one. Notice that $d_{1}\theta_{t}^{(m_l)}$ and $d_{k,1}\theta_{k,t}^{(m_l)}$ belong to the interval $[0,1)$. It is not difficult to see that $f\,:\,[0,1) \rightarrow [0,\infty)$ is monotone increasing in the interval $[0,1)$. Accordingly, it admits the continuous inverse function $g\,:\, [0,\infty)\rightarrow [0,1)$ and 
\alg{\theta_{t}^{(m_l)}&=d_1^{-1}g\left( c^{(m_l)}-o( c^{(m_l)}) \right)\nn\\
\theta_{k,t}^{(m_l)}&=d_{k,1}^{-1}g\left( c^{(m_l)}-o( c^{(m_l)}) \right).\nn} Notice that
\alg{\lim_{l\rightarrow \infty}&g\left( c^{(m_l)}-o( c^{(m_l)}) \right)=g\left(\lim_{l\rightarrow \infty} c^{(m_l)}-o( c^{(m_l)}) \right)\nn\\ &=g\left(\lim_{l\rightarrow \infty} c^{(m_l)} \lim_{l\rightarrow \infty}\left(1-\frac{o( c^{(m_l)})}{c^{(m_l)}}\right)\right)=g\left(\lim_{l\rightarrow \infty} c^{(m_l)} \right)=\lim_{l\rightarrow \infty}g\left( c^{(m_l)} \right) } 

Finally, we have 
\alg{&\lim_{l\rightarrow \infty}\frac{\|V_{t+1}^{(m_l)}\|}{\|V_{k,t+1}^{(m_l)}\|}=\lim_{l\rightarrow \infty}\sqrt{\frac{\sum_{i=1}^n \frac{1}{d_i^{-1}-\theta_{t}^{(m_l)}}}{\sum_{i=1}^n \frac{1}{d_{k,i}^{-1}-\theta_{k,t}^{(m_l)}} }}= 
\lim_{l\rightarrow \infty} \sqrt{\frac{\frac{1}{d_1^{-1}-\theta_{t}^{(m_l)}}+ \sum_{i=2}^n \frac{1}{d_i^{-1}-\theta_{t}^{(m_l)}}}{\frac{1}{d_{k,1}^{-1}-\theta_{k,t}^{(m_l)}}+  \sum_{i=2}^n \frac{1}{d_{k,i}^{-1}-\theta_{k,t}^{(m_l)}} }}\nn\\
& =\lim_{l\rightarrow \infty} \sqrt{\frac{\frac{d_1}{1-g\left( c^{(m_l)}-o( c^{(m_l)}) \right)}+ \sum_{i=2}^n \frac{1}{d_i^{-1}-\theta_{t}^{(m_l)}}}{\frac{ d_{k,1}}{1-g\left( c^{(m_l)}-o( c^{(m_l)}) \right)}+  \sum_{i=2}^n \frac{1}{d_{k,i}^{-1}-\theta_{k,t}^{(m_l)}} }}=\lim_{l\rightarrow \infty} \sqrt{\frac{\frac{d_1}{1-g\left(c^{(m_l)} \right)}+ \sum_{i=2}^n \frac{1}{d_i^{-1}-\theta_{t}^{(m_l)}}}{\frac{ d_{k,1}}{1-g\left( c^{(m_l)} \right)}+  \sum_{i=2}^n \frac{1}{d_{k,i}^{-1}-\theta_{k,t}^{(m_l)}} }}\nn\\
&=\lim_{l\rightarrow \infty} \sqrt{\frac{\frac{d_1}{1-g\left( c^{(m_l)} \right)}}{\frac{d_{k,1}}{1-g\left(c^{(m_l)} \right)} }}
} where we exploited the fact that $\lim_{x\rightarrow \infty} g(x)=1 $ in the last equality. Then, we have 
\alg{&\lim_{l\rightarrow \infty}\frac{\|V_{t+1}^{(m_l)}\|}{\|V_{k,t+1}^{(m_l)}\|}=\lim_{l\rightarrow \infty} \sqrt{\frac{\frac{d_1}{1-g\left(c^{(m_l)}\right)}}{\frac{ d_{k,1}}{1-g\left(c^{(m_l)}\right)} }}=\sqrt{\frac{ d_1}{d_{k,1}}}>0.}

Accordingly the corresponding subsequence $d_\Delta^{(m_l)}$ approaches $-\infty$ because the term $-\beta_{k,t}^{(m_l)}\| V_{t+1}^{(m_l)}\|$ dominates the logarithmic term $n\log \| V_{k,t+1}^{(m_l)}\|$. We conclude that for $c$ sufficiently large (\ref{claim_opt}) holds.
\qed

 \section{Numerical example} \label{sec_num}
 
\begin{figure}
\begin{center}
\includegraphics[width=0.4\textwidth]{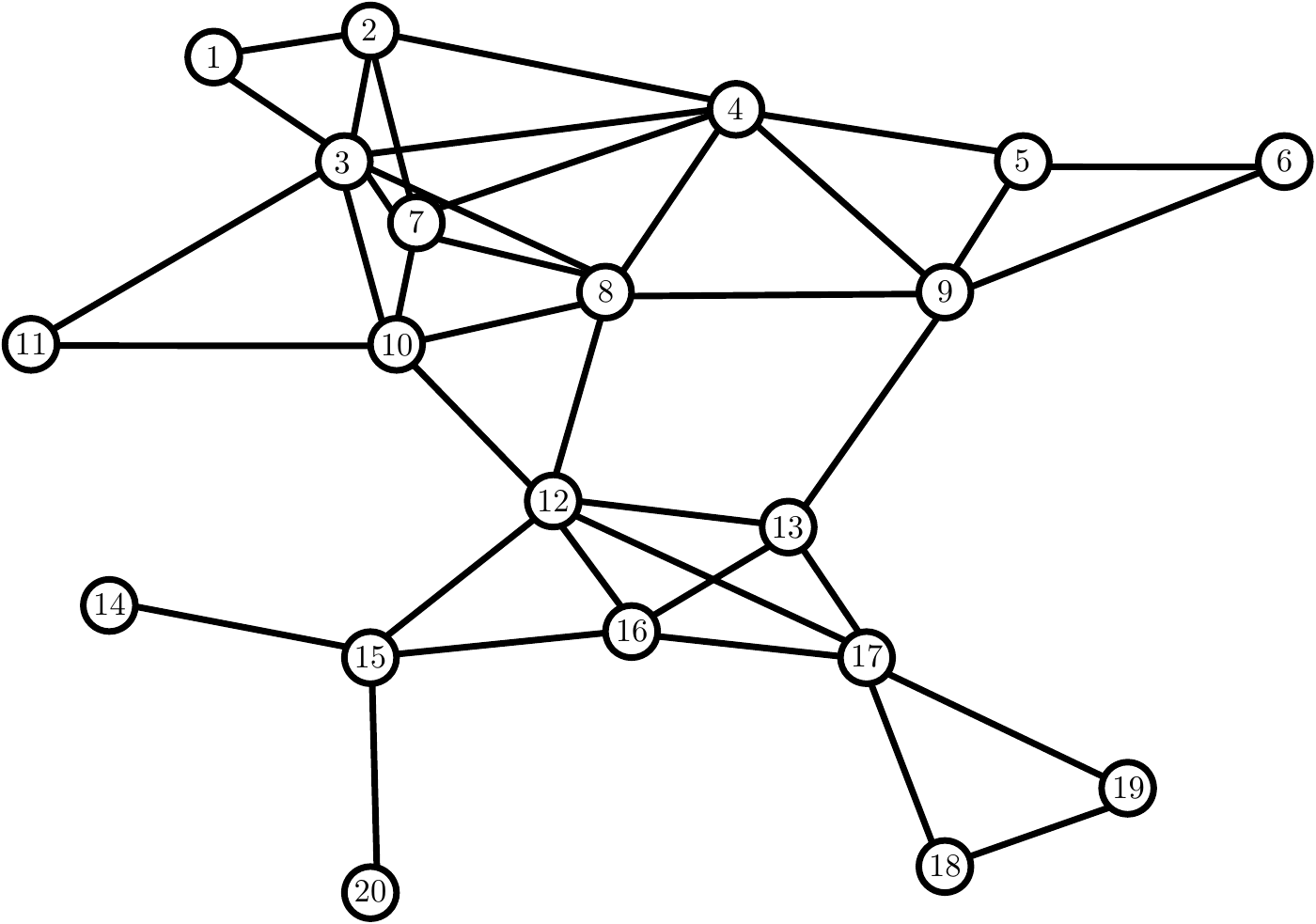}
\caption{Network of $20$ sensors for collecting noisy position measurements of the projectile.}
\label{fig_network}
\end{center}
\end{figure}

In order to evaluate the performance of the distributed robust Kalman filters, we consider the problem in \cite{cattivelli2010} of tracking the position of a projectile by using noisy position measurements obtained by a network of $N=20$ sensors depicted in Figure \ref{fig_network}.
 The model for the projectile motion is 
\alg{\label{proj_cont}\dot x^c_t =\Phi x^c_t+u^c_t}
where 
\alg{\Phi=\left[\begin{array}{cc}0 & 0 \\ I_3 & 0\end{array}\right],\nn}
$u^c_t=[\, 0 \;0 \; -g\; 0 \; 0\; 0\,]^T$, with $g=-10$, and $x^c_t=[\, v_{x,t}\;v_{y,t}\;v_{z,t}\;p_{x,t}\;p_{y,t}\;p_{z,t}\;\,]^T$ with $v$ denoting the velocity, $p$ the position and the subscripts $x,y,z$ denoting the three spatial dimensions. We discretize (\ref{proj_cont}) by using a sampling time equal to $0.1$. In this way we obtain the discrete time model $x_{t+1}=Ax_{t}+u_t$ where $x_t$ is the sampled version of $x^c_t$, $A=I_6+0.1\Phi$ and $u_t=(0.1 I_6+0.1^2\Phi/2)u^c_t$. We assume that every sensor measures the position of the projectile in either two horizontal dimensions, or a combination of one horizontal dimension and the vertical dimension (i.e. one sensor does not have measurements in all the three dimensions). Therefore, we obtain the  nominal discrete state-space model (\ref{distr_mod}) where $C_k=[\, 0\; 0 \;0 \;\mathrm{diag(1,1,0)}\,]$, in the case that the sensor measures only the horizontal positions, or $C_k=[\, 0\; 0 \;0 \;\mathrm{diag(1,0,1)}\,]$, $C_k=[\, 0\; 0 \;0 \;\mathrm{diag(0,1,1)}\,]$, in the case that the sensor measures one horizontal position and the vertical position. Moreover, we choose $B=\sqrt{0.001}I$, $R_k=D_kD_k^T=\sqrt{k} PR_0P^T$ where $R_0=0.5\cdot \mathrm{diag}(1,4,7)$ and $P$ is a permutation matrix randomly chosen for every node. Finally, the initial state $x_0$ is a Gaussian random vector with covariance matrix $P_0=I$.

In what follows, we consider the following predictors:
\begin{itemize}
\item \textsf{RKF diff} -- the distributed robust Kalman filter with diffusion step in Algorithm \ref{rob_est_diff}; the diffusion matrix $W$ is chosen as
\alg{\label{def_W_ex}w_{lk}=\left\{\begin{array}{cc}\alpha_k { n_l}, & \hbox{if } l\in\Nc_k  \\ 0, & \hbox{otherwise},\end{array}\right.}
where $n_l$ denotes the number of neighbors of node $l$ and $\alpha_k>0$ is a normalization parameter chosen in such a way that (\ref{prop_W}) holds.
\item \textsf{KF diff} -- the distributed Kalman filter with diffusion step proposed in \cite{cattivelli2010}; the diffusion matrix $W$ is chosen as in (\ref{def_W_ex}).
\item \textsf{RKF cons} -- the distributed robust Kalman filter in  Algorithm \ref{rob_est_diff} with the consensus-based update in (\ref{W_cons});  the consensus parameter is set equal to $\varepsilon=0.1$.
\item \textsf{KF cons} -- the distributed Kalman filter with consensus-based update proposed in \cite{Saber_CDC2007};  the consensus parameter is set equal to $\varepsilon=0.1$.
\item \textsf{RKF local} -- the local robust Kalman filter in Algorithm \ref{rob_est_diff} with $W=I$.
\item \textsf{KF local} -- the local Kalman filter proposed in 
\cite[p. 329]{KAILATH_OPTIMAL_FILTERING}.
\item \textsf{RKF central} -- the centralized robust Kalman filter proposed \cite{ROBUST_STATE_SPACE_LEVY_NIKOUKHAH_2013}.
\item \textsf{KF central} -- the centralized Kalman filter.
\end{itemize}

\begin{figure}
\begin{center}
\includegraphics[width=0.5\textwidth]{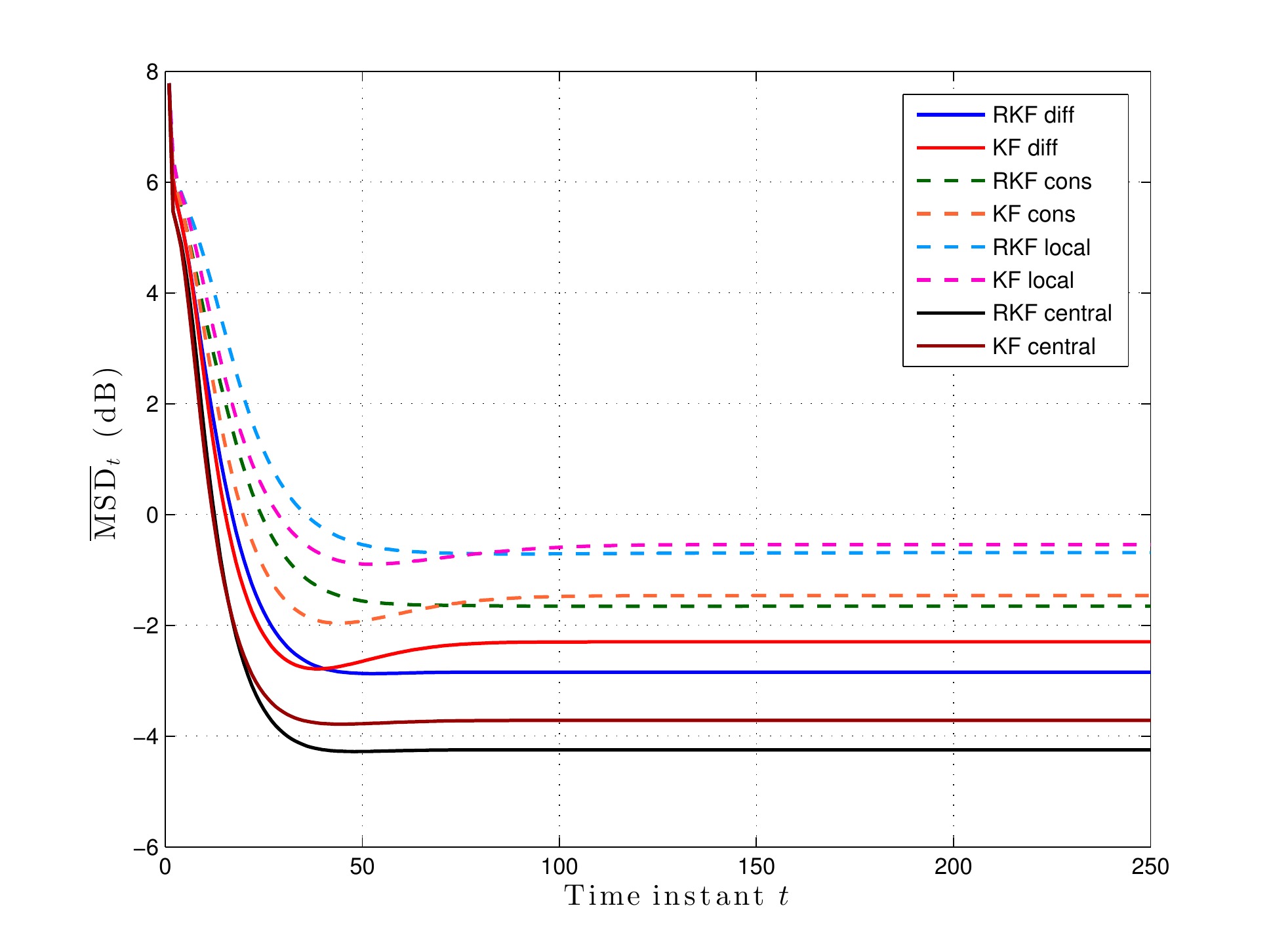}
\caption{Least favorable mean square deviation across the network with tolerance $c=0.02$.}
\label{fig_msdc1}
\end{center}
\end{figure}

\begin{figure}
\begin{center}
\includegraphics[width=0.5\textwidth]{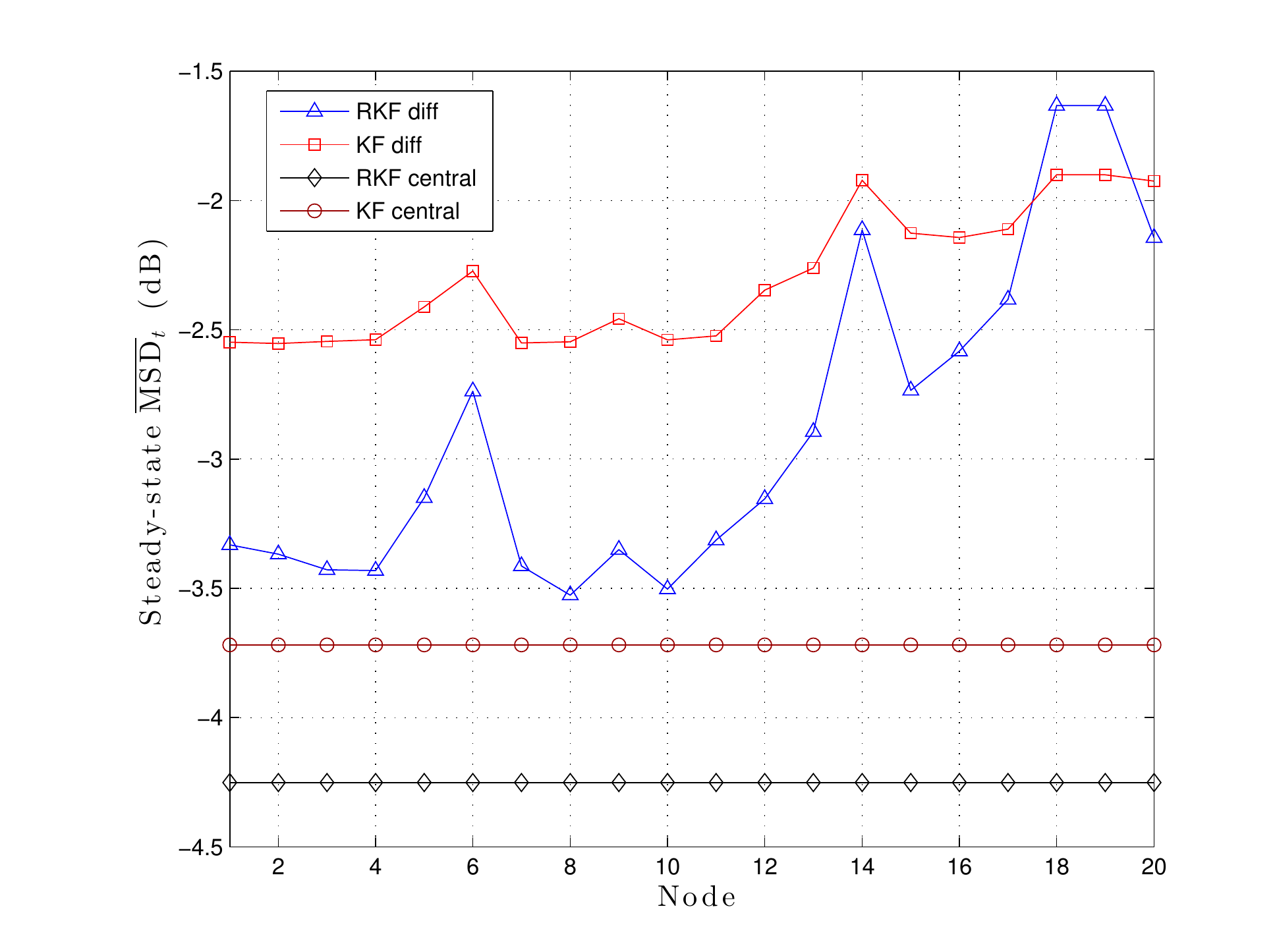}
\caption{Least favorable mean square deviation for each node in steady state with tolerance $c=0.02$.}
\label{fig_steadyc1}
\end{center}
\end{figure}

\begin{figure}
\begin{center}
\includegraphics[width=0.5\textwidth]{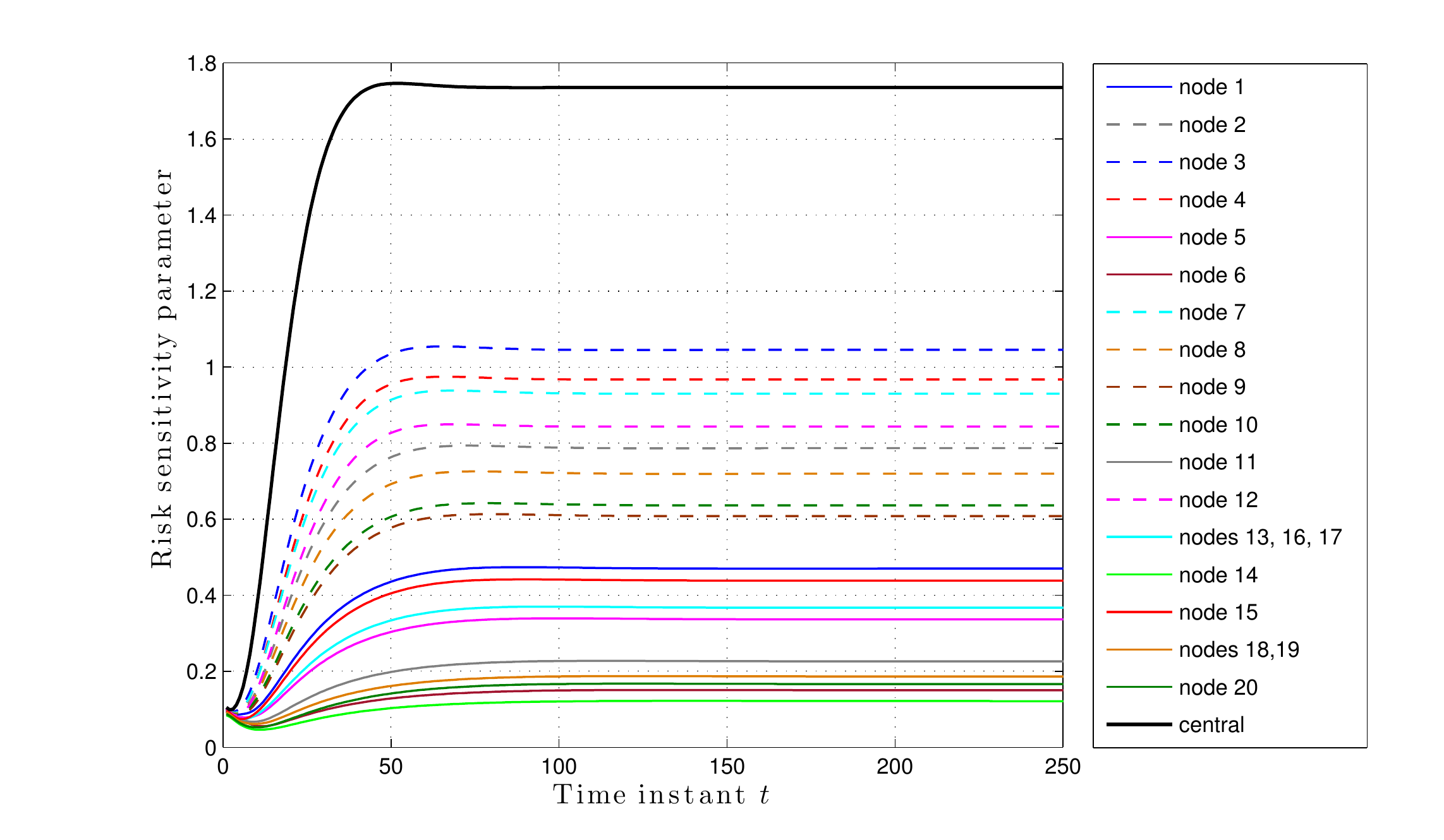}
\caption{Risk sensitivity parameters $\theta_{k,t}$, $k=1\ldots 20$, of \textsf{RKF diff} and the risk sensitivity parameter of the centralized filter \textsf{RKF} (black line) with $c=0.02$.}
\label{fig_thetac1}\end{center}
\end{figure}

In the first experiment we assume that the actual state-space model belongs to the ball defined in (\ref{ball_B}) about the aforementioned nominal model and with tolerance 
$c=0.02$. The average least favorable mean square deviation across the network is depicted in Figure \ref{fig_msdc1}.
As we can see, $\overline{\mathrm MSD}_t$ converges in steady state for any algorithm. The local algorithms \textsf{RKF local} and \textsf{KF local} provide the worst performance and the robust version behaves slightly better than the standard version in steady state. The consensus-based algorithms \textsf{RKF cons} and \textsf{KF cons} perform better than the latter and  the robust version behaves slightly better than the standard one in steady state. The diffusion-based algorithms \textsf{RKF diff} and \textsf{KF diff} provides the best distributed performance, in particular \textsf{RKF diff} performs better than \textsf{KF diff}. Finally, the centralized algorithm provides the best performance and \textsf{RKF} is the best predictor. 
The least favorable mean square deviation for each node in steady state for the diffusion-based and centralized algorithms is depicted in Figure \ref{fig_steadyc1}. 
 As we can see, \textsf{RKF diff} provides a better performance than \textsf{KF diff} in the majority of the nodes. Finally, Figure \ref{fig_thetac1} shows the risk sensitivity parameters $\theta_{k,t}$ of \textsf{RKF diff} and the risk sensitivity parameter $\theta_t$ of \textsf{RKF}.  We notice that the former are less than the latter. Therefore,  \textsf{RKF diff} reduces the risk sensitivity parameters over the network in respect to \textsf{RKF}. Such a reduction can be justified as follows. First, the larger the risk sensitivity parameter is, the more large errors are penalized, as noticed in \cite{ROBUST_STATE_SPACE_LEVY_NIKOUKHAH_2013}. Then,
it is worth observing that \textsf{RKF cons} and \textsf{RKF local} have the same risk sensitivity parameters of  \textsf{RKF diff}, indeed the value of  $\theta_{k,t}$ does not depend on the matrix $W$. So, without loss of generality, we can consider \textsf{RKF local}. \textsf{RKF local} at node $k$ and \textsf{RKF} are the same algorithm, but applied on a different state space model. 
The state space model used for  \textsf{RKF local} at node $k$ is characterized by a subset of observations of the state space model used for \textsf{RKF}. Since the mismatch modeling budget $c$ is the same for both the models, then it means that the observations of the least favorable model of \textsf{RKF} are affected by more uncertainty than the ones of \textsf{RKF local}. Accordingly, it is required to penalize large errors in \textsf{RKF} more severely than in \textsf{RKF local}, hence $\theta_t$ must be greater than $\theta_{k,t}$.

\begin{figure}
\begin{center}
\includegraphics[width=0.5\textwidth]{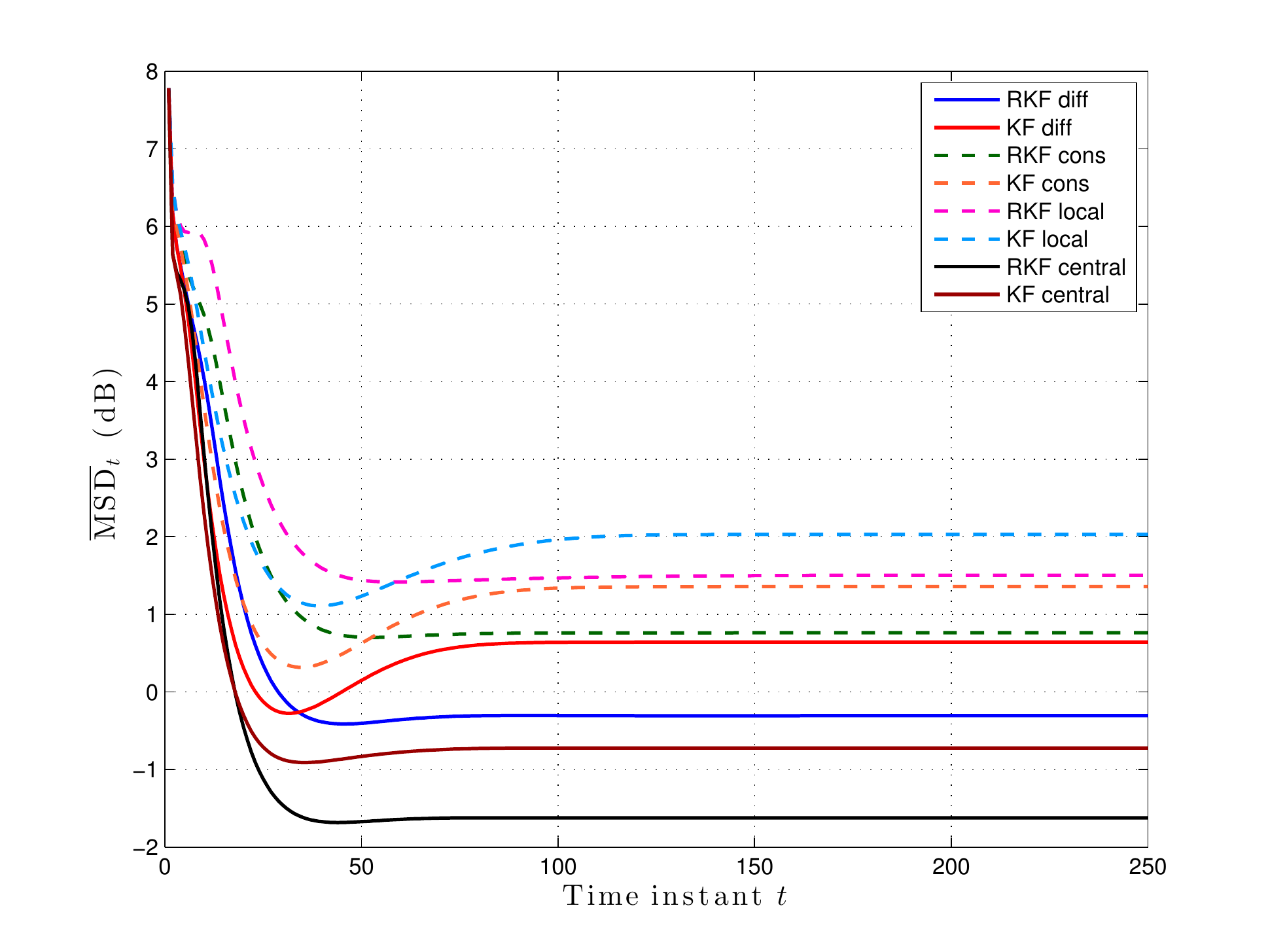}
\caption{Least favorable mean square deviation across the network with tolerance $c=0.06$.}
\label{fig_msdc2}\end{center}
\end{figure}

\begin{figure}
\begin{center}
\includegraphics[width=0.5\textwidth]{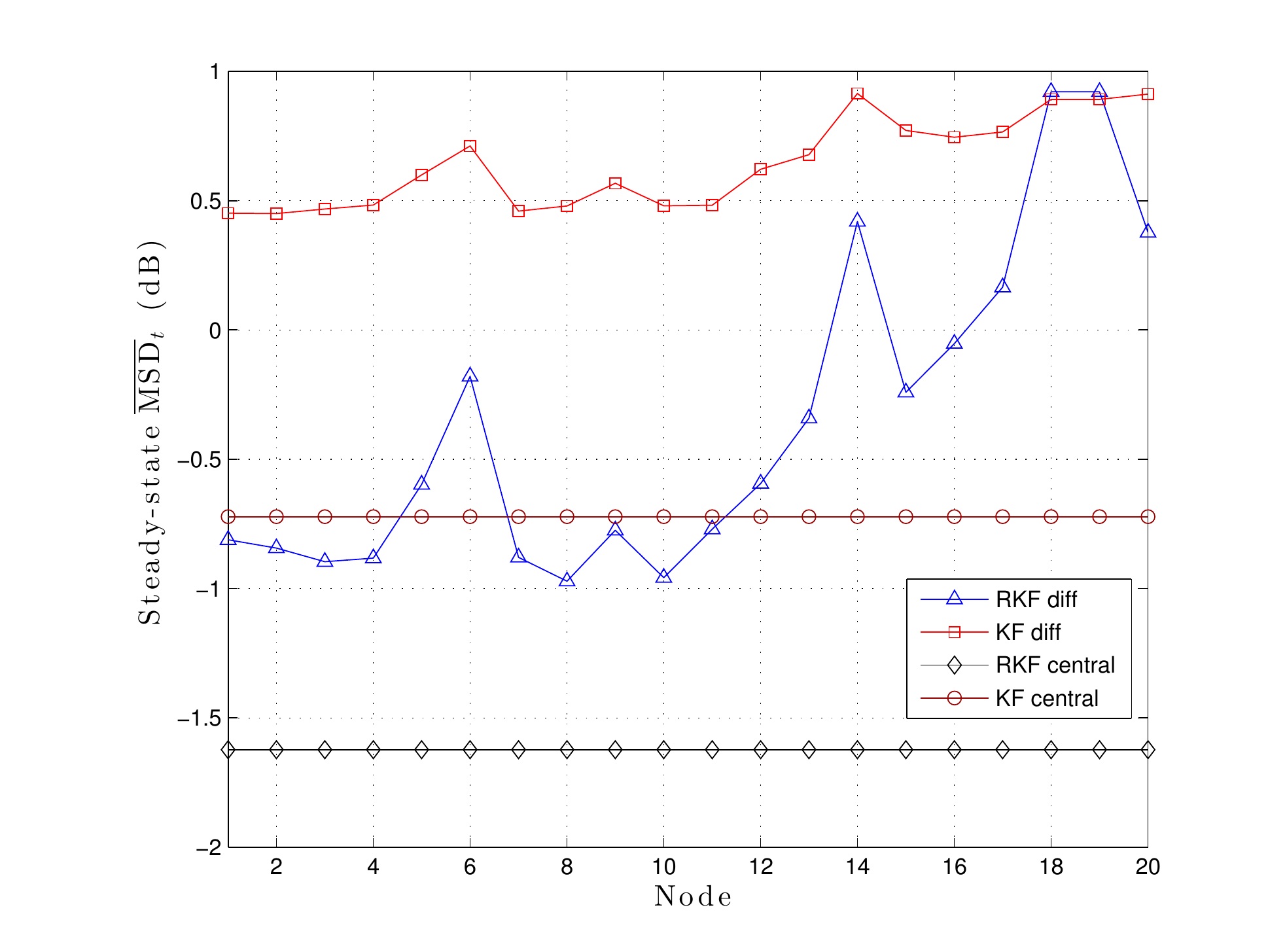}
\caption{Least favorable mean square deviation for each node in steady state with tolerance $c=0.06$.}
\label{fig_steadyc2}\end{center}
\end{figure}

\begin{figure}
\begin{center}
\includegraphics[width=0.5\textwidth]{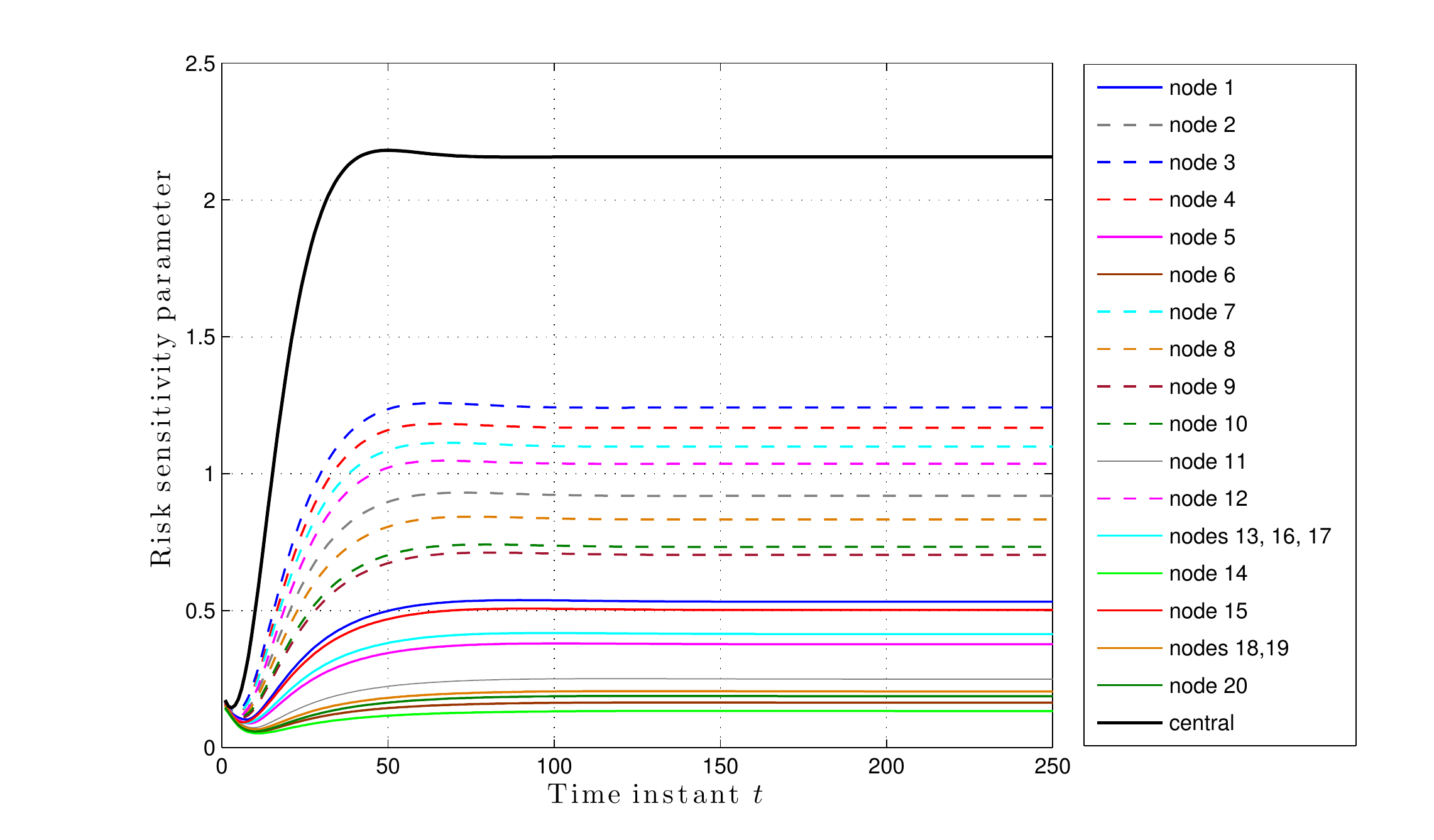}
\caption{Risk sensitivity parameters $\theta_{k,t}$, $k=1\ldots 20$, of \textsf{RKF diff} and the risk sensitivity parameter of the centralized filter \textsf{RKF} (black line) with $c=0.06$.}
\label{fig_thetac2}\end{center}
\end{figure}

In the second experiment we have considered the case that in the actual model there are large deviations in respect to the nominal one. More precisely, we have chosen $c=0.06$. The least favorable mean square deviation across the network is depicted in Figure \ref{fig_msdc2}. As we can see, all the robust distributed algorithms outperform the corresponding standard distributed algorithms in steady state. Among the robust algorithms, \textsf{RKF diff} gives the best performance, then we have \textsf{RKF cons} and 
finally \textsf{RKF local}. The least favorable mean square deviation for each node in steady state for the diffusion-based and centralized algorithms is depicted in Figure \ref{fig_steadyc2}. 
 As we can see, \textsf{RKF diff} provides a better performance than \textsf{KF diff} in almost all nodes: the unique exception regards  two nodes wherein  \textsf{RKF diff} performs slightly worse than \textsf{KF diff}. The risk sensitivity parameters $\theta_{k,t}$ of \textsf{RKF diff} and the risk sensitivity parameter $\theta_t$ of \textsf{RKF} are depicted in Figure \ref{fig_thetac2}.   Also in this case, \textsf{RKF diff} reduces the risk sensitivity parameters over the network in respect to \textsf{RKF}. On the hand, the values of all these risk sensitivity parameters has been increased in respect to the case $c=0.02$. Indeed, in the current case the mismatch modeling budget has been increased and thus it is required to penalize large error more severely.

\section{Conclusions} \label{sec_concl}
In this paper, we have considered a filtering problem over a sensor network and under model uncertainty. We have proposed a robust distributed algorithm with diffusion step. We have derived the least favorable performance for this algorithm and showed that the least favorable mean square deviation across the network does converge to a finite constant value provided that the mismatch modeling budget allowed for each time step is sufficiently small. Finally, a numerical example showed that this robust algorithm is preferable than the standard one in the presence of model deviations.

\end{document}